\documentclass{amsart}

\usepackage[all]{xy}
\usepackage{amsmath,amsthm,amssymb,amscd,enumitem,mathtools}
\usepackage{amsfonts}
\usepackage{rotating}
\usepackage{euscript}
\usepackage{pst-node}
\usepackage{epsfig,verbatim}
\def\be{\begin{equation}}
\def\ee{\end{equation}}

\def\C{{\mathbb C}} 
\def\f{\EuScript}
\def\N{{\mathbb N}} 
\def\P{{\mathbb P}}

\def\R{{\mathbb R}} 
\def\Q{{\mathbb Q}}

\def\phi{{\varphi}}
\def\v{{\varepsilon}} 

\def\deg{{\rm deg\,}} 

\def\Aut{{\rm Aut}}

\def\GCD{{\rm GCD }}
\def\LCM{{\rm LCM }}

\def\bp{\begin{proposition}}
\def\ep{\end{proposition}}

\def\bt{\begin{theorem}}
\def\et{\end{theorem}}
\def\br{\begin{remark}}
\def\er{\end{remark}}
\def\be{\begin{equation}}
\def\bee{\begin{equation*}}
\def\l{\label}

\def\m{\mu}
\def\ee{\end{equation}}
\def\eee{\end{equation*}}
\def\bl{\begin{lemma}}
\def\el{\end{lemma}}
\def\bc{\begin{corollary}}
\def\ec{\end{corollary}}
\def\pr{\noindent{\it Proof. }}

\def\bd{\begin{definition}}
\def\ed{\end{definition}}

\mathtoolsset{showonlyrefs}

\newtheorem{theorem}{Theorem}[section]
\newtheorem{lemma}[theorem]{Lemma}
\newtheorem{definition}[theorem]{Definition}
\newtheorem{corollary}[theorem]{Corollary}
\newtheorem{proposition}[theorem]{Proposition}

\theoremstyle{definition}

\theoremstyle{definition}
\newtheorem{remark}[theorem]{Remark}

\begin{document}
\title{On amenable semigroups of rational functions}
\author{Fedor Pakovich}
\thanks{
This research was supported by ISF Grant No. 1432/18}
\address{Department of Mathematics, Ben Gurion University of the Negev, Israel}
\email{
pakovich@math.bgu.ac.il}

\begin{abstract} 
We characterize left and right amenable semigroups of polynomials of one complex variable with respect to the composition operation.
We also prove a number of results about amenable semigroups of arbitrary rational functions. In particular, we 
show that under quite general conditions a semigroup of rational functions is amenable if and only if it is a subsemigroup 
of the centralizer of some  rational function. 
\end{abstract}

\maketitle

\section{Introduction} 

The concept of  amenable  group was introduced by von Neumann in 1929 in the paper \cite{von}.   Defined initially in terms of invariant measures in relation with the Banach-Tarski paradox, nowadays the group amenability is known to be equivalent to variety of different conditions 
and to have connections to numerous branches of mathematics (see e.g. \cite{wagon}, \cite{ju} and the bibliography therein). 
The notion of amenability  was extended to 
semigroups by Day \cite{day}, who also introduced the term itself. Naturally, the absence of inverse elements in semigroups requires substantial changes in definitions and leads to new phenomenons. For example, a 
semigroup can be left amenable but not right amenable, 
amenable semigroups can contain non-amenable subsemigroups etc. (see e.g. \cite{pater}).

Let us recall that a semigroup $S$ is called {\it left amenable} if it admits a finitely additive probability measure $\mu$ defined on all the subsets of $S$  
such that for all $a\in S$ and $T\subseteq S$ 
the equality \be \l{perd0} \mu(a^{-1}T) = \mu(T)\ee holds, where   
 the set $a^{-1} T$ is defined by the formula  $$a^{-1} T=\{s \in S\, | \,as \in T\}.$$ 
Equivalently, $S$ is left amenable if there is a mean on $l_{\infty}(S)$, which is invariant under the natural
left action of $S$ on the dual space $l_{\infty}(S)^*$ (see e.g. \cite{pater}). The right amenability is defined similarly. A semigroup is called {\it amenable} if 
there exists a mean on $l_{\infty}(S)$, which is invariant under the left and the right action of $S$ on $l_{\infty}(S)^*$. By the theorem of Day (see \cite{day0}, \cite{day}), this is equivalent to the condition that $S$ is left and right amenable. 

In this paper, we investigate the amenability of semigroups of polynomials and more generally of rational functions of one complex variable  with respect to the composition operation. 
To our best knowledge, for the first time this topic was investigated only recently  by Cabrera and Makienko in the  paper \cite{peter}. Among other things, they proved 
that if $S$ is a semigroup generated by finitely many polynomials at least one of which is of degree greater than one and is not  conjugate to $z^n$ or $\pm T_n$, then 
$S$ is right amenable if and only if there exists an  $S$-invariant probability measure $\mu$ on $\C\P^1$ such that the measure of maximal entropy of every  element of $S$ of degree at least two coincides with $\mu.$  Cabrera and Makienko also provided conditions for amenability of  semigroups of polynomials and more generally of semigroups of rational functions if, an addition to the amenability, some extra conditions are satisfied.

In this paper, 
we obtain a full description of left and right amenable semigroups of polynomials, complementing  and generalizing the results of \cite{peter}. In particular, we  show that for  semigroups of polynomials the left 
amenability implies the right amenability, and that amenability conditions are equivalent to weaker algebraic  {\it reversibility} conditions. To formulate our results  explicitly we introduce several definitions. 
Let us recall that a semigroup $S$ is called {\it left reversible}  if for all $a,b\in S$ the right ideals 
$aS$ and  $bS$ have a non-empty intersection, that is, if for all $a,b\in S$  there exist $x,y\in S$ such that \be \l{reva} ax=by.\ee
It is well-known and follows easily from the definition  that any left amenable semigroup is left reversible.

For a  rational function $R$ of degree at least two, we denote by $C(R)$ the semigroup of  rational functions commuting with $R$, and by $ G(R)$ the group of M\"obius transformations $\sigma$ such that $R\circ \sigma=\nu \circ R$ for some  M\"obius transformations $\nu$. 
 It is easy to see that  $ G(R)$ is indeed a group 
and that the map $\gamma_R:\sigma \rightarrow \nu_{\sigma}$ is a homomorphism from   $ G(R)$ to the group $\rm{Aut}(\C\R^1)$. 
For  a subgroup $\Gamma$  of $G(R)$ such that $\gamma_R(\Gamma)\subseteq \Gamma$, we denote by  $S_{\Gamma,R}$ the semigroup of rational functions generated by $\Gamma$ and $R$. 
 We will say that 
a polynomial $P$  is {\it special} if it is conjugate to $z^n,$ $n\geq 2,$ or $\pm T_n,$ $n\geq 2,$ where $T_n$ is the $n$th Chebyshev polynomial.

 In this notation, our main result about left amenable and amenable semigroups of polynomials is  following.

\bt \l{-1} Let $S$ be a semigroup of  polynomials containing at least one non-special polynomial of degree greater than one. Then  the following conditions are equivalent:

\begin{enumerate} [label=\arabic*)]

\item The semigroup $S$ is left reversible.

\item The semigroup $S$ is left amenable. 

\item The semigroup $S$ is amenable.

\item The semigroup $S$ is a subsemigroup of $S_{\Gamma,R}$ for some non-special polynomial $R$  
of degree at least two  and a group $\Gamma\subseteq G(R)$ such that $\gamma_R(\Gamma)= \Gamma$.

\item The semigroup $S$ is a subsemigroup of  $C(P)$ for some  non-special polynomial $P$  
of degree at least two.

\end{enumerate}

\noindent Furthermore, if $S$ contains only polynomials of degree at least two, then any of the above conditions is equivalent to the 
condition that  for all $A,B\in S$ there exist $k,l\geq 1$ such that $A^{\circ k}=B^{\circ l}$.

\et

To formulate the analogue of Theorem \ref{-1} for right amenable semigroups of polynomials we introduce two other types of special semigroups. 
Let us recall that by the results of Freire, Lopes, Ma\~n\'e (\cite{flm}) and Lyubich (\cite{l}), 
 for every rational function $P$ of degree $n\geq 2$ there exists a unique probability measure $\mu_P$ on $\C\P^1$, which is invariant under $P$, has support equal to the Julia set $J_P$, and achieves maximal entropy $\log n$ among all $P$-invariant probability measures.

For a  rational function $P$ of degree at least two,  we denote by $E(P)$ the semigroup consisting of  rational functions $Q$ of degree at least two such that $\mu_Q=\mu_P$, completed by  $\mu_P$-invariant M\"obius transformations. Finally, for a compact set $K\subset \C$, 
we denote by $I(K)$ the semigroup of all polynomials $A$ satisfying $A^{-1}\{K\}=K.$  

Our main result about right amenable semigroups of polynomials is  following.

\bt \l{-11}  Let $S$ be a semigroup of  polynomials containing at least one non-special polynomial of degree greater than one. Then  the following conditions are equivalent:

\begin{enumerate} [label=\arabic*)]

\item The semigroup $S$ is right reversible.

\item The semigroup $S$ is right amenable.

\item The semigroup $S$ is a subsemigroup of $S_{\Gamma,R}$ for some non-special polynomial $R$  
of degree at least two  and a group $\Gamma\subseteq G(R)$ such that $\gamma_R(\Gamma)\subseteq \Gamma$.

\item The semigroup $S$ is subsemigroup of  $I(K)$ for some compact set $K\subset \C$, which is neither a union of concentric circles nor a segment.

\item The semigroup $S$ is a subsemigroup of  $E(P)$  for some non-special polynomial $P$  
of degree at least two.

\item The semigroup $S$ contains no free subsemigroup of rank two.

\end{enumerate} 

\noindent Furthermore, if $S$ contains only polynomials of degree at least two, then any of the above conditions is equivalent to the condition that  for all $A,B\in S$ there exist $k,l\geq 1$  such that $A^{\circ 2k}=A^{\circ k}\circ B^{\circ l}$ and $B^{\circ 2l}=B^{\circ l}\circ A^{\circ k}$. 
\et

Notice that  for semigroups of polynomials of degree at least two the fifth condition of Theorem \ref{-11} is equivalent to the requirement that all elements of $S$ share a measure of maximal entropy. In particular, 
Theorem \ref{-11} confirms in the polynomial case the following conjecture proposed  in  \cite{peter}: if a semigroup of rational functions  of degree at least two $S$ contains no free subsemigroup of rank two, then all elements of $S$ share a  measure of maximal entropy.  Theorem \ref{-11} also implies the following statement in spirit of von Neumann conjecture for amenable groups:  if a semigroup
of polynomials $S$ is not right amenable, then $S$  has a free subsemigroup of rank two.
Finally, since in the polynomial case having the same measure of maximal
entropy is equivalent to having the same Julia set, Theorem \ref{-11} implies that if the Julia sets of polynomials $A$ and $B$ are not equal, then the semigroup $\langle A,B \rangle$ contains a free subsemigroup of rank two.

In addition to the polynomial case, we study the amenability of semigroups of  arbitrary rational functions and prove a partial generalization of Theorem \ref{-1} to a wide class of such functions. 
Recall that a semigroup $S$ is called {\it left cancellative} if the equality $ab=ac,$ where $a,b,c\in S,$ implies that $b=c$. Right cancellative semigroups are defined similarly. Any semigroup of rational functions is obviously right cancellative but not necessarily left cancellative. 
Following \cite{tame}, we say that a rational function $A$ of degree at least two is {\it tame} if the algebraic curve $$A(x)-A(y)=0$$ has no factors of genus zero or one distinct from the diagonal. By the Picard theorem, this condition  is equivalent to the condition that the equality   
$$ A\circ f=A\circ g,$$ where  $f$ and $g$  are functions  meromorphic on $\C$, implies that $f\equiv g.$
We say that a semigroup of rational functions $S$ is tame if $S$ consists of tame rational functions. 
Clearly, any tame semigroup of rational functions $S$ is cancellative, so the  tameness condition can be regarded as 
a  strengthening of the cancellativity condition. Notice that tame rational functions form a subsemigroup of $\C(z).$ 

In the above notation, our main result about left amenable
semigroups of rational functions is following.  

\bt \l{11} Let $S$ be a tame semigroup of rational functions. Then the following conditions are equivalent.

\begin{enumerate} [label=\arabic*)]

 \item The semigroup $S$   is    left reversible.

 \item The semigroup $S$   is   left  amenable.

 \item The semigroup $S$   is   amenable. 

 \item The semigroup $S$   is  a subsemigroup of $C(P)$  
 for some  tame rational \linebreak function $P$.

\item For all $A,B\in S$ there exist $k,l\geq 1$ such that $A^{\circ k}= B^{\circ l}$.

 \item The semigroup $S$ contains no free subsemigroup of rank two.


\end{enumerate} 

\et 

\enlargethispage*{0.15cm}

Our approach to the study of left amenable semigroups of rational functions relies on using the reversibility condition. Specifically, applying condition \eqref{reva} to powers of $a$ and $b$, we conclude that 
if a semigroup of rational functions $S$ is left amenable, then for all $A,B\in S$ of degree at least two  the algebraic curves 
\be \l{cu1} A^{\circ n}(x)-B(y)=0, \ \ \ n\geq 1,\ee
and, more generally,  the algebraic curves 
\be \l{cu2} A^{\circ n}(x)-B^{\circ m}(y)=0, \ \ \ n\geq 1, \ \ \ m\geq 1,\ee have a factor of genus zero.

The problems of describing pairs of rational functions satisfying the above conditions arose 
recently in the context of arithmetic dynamics. 
Specifically, the problem of describing  $A$ and $B$ such that all curves \eqref{cu1} have a factor of genus zero or one 
is a geometric counterpart of the following problem of the arithmetic nature posed in \cite{jones}:   which rational functions $A$
defined over a number field $K$ have a $K$-orbit containing infinitely many points from the value set $B(\P^1(K))$? These problems have been studied in \cite{jones}, \cite{hyde}, \cite{aol}.  In particular, in \cite{aol}, a description of such $A$ and $B$  in terms of semiconjugacies and Galois coverings was obtained. 

In turn, the problem of describing pairs of rational functions  $A$ and $B$ such that all curves \eqref{cu2} have a factor of genus zero or one is a geometric counterpart of the problem of describing pairs of rational functions $A$ and $B$ having orbits
with infinite intersection. In case $A$ and $B$ are polynomials, the last problem was solved in the papers \cite{gtz}, \cite{gtz2}, where it was shown that such orbits
 exist if and only if $A$ and $B$ have a common iterate.   
This result was extended to tame rational functions in the paper \cite{tame}, and our approach to the proof of Theorem \ref{11} is based on ideas and results of this paper.

Notice that  in the context of right amenability  the analogues of the above problems about algebraic curves
can be formulated in terms of intersections of subfields of rational functions as follows: given rational functions $A$ and $B$, under what conditions the fields   $$\C(A^{\circ n})\cap \C(B), \ \ \ n\geq 1,$$ and, more generally, the fields  
$$\C(A^{\circ n})\cap \C(B^{\circ m}), \ \ \ n\geq 1,$$ contain a non-constant rational function?
These problems however have a different flavor, and are not considered in this paper.

The last class of semigroups of rational functions whose amenability is characterized in this paper is the class of semigroups  $S_{\Gamma,R}$. 
In a sense, these semigroups  are the simplest examples of non-cyclic semigroups of rational functions, and the polynomial case suggests that they  play an important role in the whole theory. Our main result concerning semigroups  $S_{\Gamma,R}$ is following.

\bt \l{pi} Let $R$ be a  rational function of degree $n\geq 2$  not conjugate to $z^{\pm n}$,  and $\Gamma$  a subgroup of $G(R)$ such that 
$\gamma_R(\Gamma)\subseteq  \Gamma$. Then every subsemigroup of $S_{\Gamma,R}$ is  right amenable. On the other hand, 
$S_{\Gamma,R}$ is left amenable if and only if $\gamma_R(\Gamma)= \Gamma$. Moreover, in the last case every subsemigroup of $S_{\Gamma,R}$ is amenable. 
\et

The paper is organized as follows. In the second section, after recalling some basic definitions and results about abstract amenable semigroups, we study
 semigroups $C_{\infty}(P)$ consisting   of rational functions commuting with some {\it iterate} of  a fixed rational function $P$ of degree at least two. Assuming that $P$ is not a Latt\`es map and is not conjugate to  $z^{\pm n}$  or $\pm T_n,$ we describe basic properties of such semigroups basing on results about commuting rational functions  from the papers \cite{r} and \cite{rev}. In  particular, we prove the amenability of $C_{\infty}(P)$ and all its subsemigroups. 

In the third section, we briefly discuss properties of semigroups $E(P)$ basing on the link between rational functions sharing a measure of maximal entropy 
and the system of functional equations 
$$A^{\circ 2k}=A^{\circ k}\circ B^{\circ l},\ \ \ \ \ \ \ \  B^{\circ 2l}=B^{\circ l}\circ A^{\circ l},$$  established by 
Levin and Przytycki (\cite{lev}, \cite{lp}). 
In the fourth section, we study the semigroups $S_{\Gamma,R}$, and prove Theorem \ref{pi}. 
We also show that
$S_{\Gamma,R}\subseteq E(R)$, and that  $S_{\Gamma,R}\subseteq C_{\infty}(R)$, whenever $\gamma_R(\Gamma)= \Gamma$.  In the fifth section, we study tame semigroups of rational functions and prove an extended version of Theorem \ref{11}. The proof is based on results of the paper \cite{tame} completed by the following stabilization result of independent interest: for a tame rational function $P$ the semigroup $C_{\infty}(P)$ coincides with the semigroup $C(P^{\circ k})$ for some $k\geq 1$.

Finally, in the sixth and the seventh sections,  we consider semigroups of polynomials. Specifically,  in the sixth section, we 
characterize reversible and Archimedean semigroups of polynomials.  Then, in the seventh section, using the results of the papers \cite{p1} and \cite{pj}, we characterize the semigroups $C(P)$ and $I(K)$. 
Finally, we prove extended versions of Theo\-rem \ref{-1} and Theorem \ref{-11}.

\section{Amenability of semigroups $C_{\infty}(f)$ and their subsemigroups} 
\subsection{Amenable semigroups}\l{xcv} 
We recall that a semigroup $S$ is called {\it left amenable} if it admits a finitely additive probability
measure $\mu$, defined on all the subsets of $S$, which is left invariant in the following sense. For all $T\subseteq S$ and $a \in S$ the equality 
\be \l{perd} \mu(a^{-1}T) = \mu(T)\ee holds, where the set $a^{-1} T$ is defined by the formula  $$a^{-1} T=\{s \in S\, | \,as \in T\}.$$ 
Equivalently, $S$ is left amenable if there is a mean on $l_{\infty}(S)$ which is invariant under the natural
left action of $S$ on the dual space $l_{\infty}(S)^*$ (see e.g. \cite{pater}).  The right amenability is defined similarly. A semigroup is called {\it amenable} if 
there exists a mean on $l_{\infty}(S)$, which is invariant under the left and the right action of $S$ on $l_{\infty}(S)^*$. By the theorem of Day (see \cite{day0}, \cite{day}), this is equivalent to the condition that $S$ is left and right amenable, and in this paper we will use the last condition as the definition of amenability.  For a given semigroup $S$, the left or the right amenability does not imply in general the opposite type of amenability. However,   any statement concerning the left amenability of semigroups has a ``right'' counterpart, which is obtained by   
switching between a semigroup $S$ with a binary operation $f(x,y)$ and a  semigroup $S'$ with the same set of elements and a binary operation $f'(x,y)=f(y,x)$.  

We start by recalling some definitions and results concerning abstract amenable semigroups. Mostly, we will discuss the ``left'' case, leaving the formulations 
in the ``right'' case to the reader. Nevertheless, the results used only in the ``right'' case will be given accordingly. 

The following statement lists some types of amenable and not amenable semigroups (see \cite{day}, \cite{pater}).

\bt \l{xrukas} Every abelian semigroup is amenable. Every finite group is ame\-nable. The free semigroup of rank two is not left or right amenable.
\qed 
\et

We recall that a semigroup $S$ is called  {\it left  cancellative} if the equality $$ab=ac$$  for $a,b,c\in S$ implies the equality $b=c$. A semigroup $S$ is said to satisfy the {\it left F{\o}lner condition} if for every finite subset $H$  of 
 $S$ and every $\varepsilon >0$ there is a finite subset $F$ of $S$ with $$\vert sF\setminus F\vert \leq \varepsilon \vert F\vert$$ for all $s\in H.$ If for every finite subset $H$  of 
 $S$ and every $\varepsilon >0$ there is a finite subset $F$ of $S$ with 
  $$\vert F\setminus sF\vert \leq \varepsilon \vert F\vert$$ for all $s\in H,$ then  $S$ is said to satisfy the {\it strong  left F{\o}lner condition}. It is known that the   strong F{\o}lner condition implies the left amenability (\cite{aw}), while the left amenability implies the F{\o}lner condition (\cite{f}, \cite{n}).
In case  $S$ is left cancellative, the sets $F$ and $sF$ have the same cardinality, implying that 
$$\vert sF\setminus F\vert=\vert F\setminus sF\vert.$$ Thus, the following criterion holds. 

\bt \l{fo} A left cancellative semigroup is left amenable if and only if it satisfies the left F{\o}lner condition. \qed 
\et

In addition to Theorem \ref{fo}, we will use the following criterion (see \cite{day}, p. 516).

\bl \l{da} 
If $\Sigma_n$ is a set of left amenable subsemigroups in a semigroup $\Sigma$ such that for every $m,n$ there exists $p$ such that $\Sigma_n,\Sigma_m\subseteq \Sigma_p$ and $\Sigma=\bigcup_{i=1}^{\infty} \Sigma_i,$ then $\Sigma$ is left amenable. 
\el

We recall that a semigroup $S$ is called {\it left reversible}  if for all $a, b\in S$ the condition $$aS\cap bS \neq \emptyset$$ holds, or, equivalently, if for all $a,b\in S$ there exist $x,y\in S$ such that \be \l{suk} ax=by.\ee

The following statement is obtained easily  from the definitions  (see \cite{pater}, Proposition 1.23).

\bp \l{pater}
Every left amenable semigroup is left reversible. \qed
\ep

In distinction with the group case, a subsemigroup of a left amenable
semigroup or even of an amenable group is not necessarily left amenable.
 However, the following result holds (see \cite{f}, \cite{don}). 

\bt \l{fr} Let $S$ be a cancellative semigroup such that $S$ contains no free
subsemigroup on two generators. If $S$ is left amenable, then every subsemigroup
of $S$ is left amenable. \qed 
\et 
   
For a semigroup $U$, we denote by  ${\rm End}(U)$ the set of endomorphisms of $U$. 
Suppose that $U$ and $T$ are semigroups with a homomorphism
$\rho: T \rightarrow {\rm End}(U).$  Denoting for $a\in T$ the endomorphism $\rho(a)$ of $U$  by $\rho_a,$  we define the semidirect product of $U$ by $T$
as the semigroup $S=U\underset{\rho}{\times} T$ of ordered pairs $(u, a)$, where 
$u\in U$ and $a\in T$, with the operation 
$$(u, a)(v, b) = (u\rho_a(v), ab).$$ 

The following result was proved in \cite{kla}. 

\bt \l{kla} If $U$ and $T$ are right amenable semigroups
with a homomorphism $\rho: T \rightarrow {\rm End}(U)$, then $S=U\underset{\rho}{\times} T$ is right
amenable. 
\qed 
\et 
  
Let us recall that a congruence on a semigroup $S$ is an equivalence relation on $S$ compatible with the structure of semigroup. 
Let $S$ be a right reversible semigroup, and let $\sim$ be the relation on $S$, which identifies $x$ and $y$ if there exists $s\in S$ for which \be \l{listik} s\circ x=s\circ y.\ee In this notation, the following criterion for the right amenability holds (see  \cite{pater}, Proposition 1.24 and Proposition 1.25).

\bt \l{lastik} Let $S$ be a right reversible semigroup. 
Then the relation $\sim$ is a congruence on $S$ and the semigroup $S/\sim$ is left cancellative. Moreover, $S$  is right amenable if and only if $S/\sim$ is right amenable. \qed 
\et 

Finally, we need the following simple statement. 

\bl \l{xrun}  Let $S$ be a left cancellative semigroup that contains no free subsemigroup of rank two. Then $S$ is left reversible. 
\el 
\pr Let $a\neq b$ be elements of $S$. By condition, the semigroup generated by $a$ and $b$ is not free. Therefore, there exist $x_1,x_2,\dots ,x_k\in \{a,b\}$ and   
$y_1,y_2,\dots ,y_l\in \{a,b\}$ such that \be \l{sok} x_1x_2\dots x_k=y_1y_2\dots y_l,\ee but the words in $a$ and $b$ in the parts of this equality are different. 
It follows  from the left cancellativity that without loss of generality we may assume that $x_1\neq y_1$. Thus, either $x_1=a$ and $y_1=b$, or $x_1=b$ and $y_1=a.$ Moreover, the condition $a\neq b$ implies that $k\geq 2$ and $l\geq 2$. Thus, if, say, $x_1=a$ and $y_1=b$, the elements  
$$x= x_2\dots x_k, \ \ \ \ \   y=y_1y_2\dots y_l $$ of $S$ satisfy \eqref{suk}. \qed

\subsection{Archimedean, power joined, and power twisted semigroups of rational functions} 
Let $F$ and $G$ be rational functions. 
 We say that $G$ is a {\it compositional left factor} of a $F$, if $F = G\circ H$ for some rational function $H$. Compositional right factors are defined in a similar way. 
For a semigroup of rational functions $S$, we denote by $\underline{S}$ and by $\overline{S}$  
the subsets of $S$ consisting of rational functions of degree one and of degree greater than one, correspondingly. It is easy to see that $\underline{S}$ and $\overline{S}$  are subsemigroups of $S$.

We recall that a semigroup $S$ is called  {\it power joined} if for all $a,b\in S$ there exist $k,l\geq 1$ such that 
\be \l{eba0} a^k=b^l,\ee 
and it is called  {\it power twisted} if for all $a,b\in S$  there exist $k_1,k_2,l\geq 1$  such that 
\be \l{pizz} a^{ k_1}=a^{ k_2} b^{l}.\ee Since \eqref{eba0} implies that  $a^{k+1}=ab^l,$ any power joined semigroup is power twisted, but the inverse is not true in general. 
A semigroup $S$ is called left (resp. right) {\it Archimedean} if for all $a,b\in S$ there exist $n\geq 1$ and $u\in S$ such that $ a^n=ub$ (resp. $a^n=bu$).

\bl \l{svin1} Let $S$ be a semigroup. If $S$  is  power joined, then $S$ is left and right  Archimedean. On the other hand, if $S$  is power twisted, then $S$ is left Archimedean. Finally, 
if $S$ is right (resp. left) Archimedean, then $S$ is left (resp. right) reversible.
\el
\pr By definition, if $S$  is  power joined, then for all $a,b\in S$ there exist $k,l\geq 1$ such that  \eqref{eba0} holds. Moreover, since equality \eqref{eba0} implies the equality $a^{ 2k}=b^{ 2l}$, 
without loss of generality we may assume that  $l\geq 2$, implying that the equalities  \be \l{eq1} a^{{k}}=bu, \ \ \ \   a^{ {k}}=ub \ee 
hold for  $ u=b^{ (l-1)}\in S$  (the assumption $l\geq 2$ is necessary since $S$ may not  contain a unit element, and thus the expression $b^{l-1}$ is defined  only for $l>1$). 

Similarly, the second part of the lemma is true since equality \eqref{pizz} implies 
the equality
$$  a^{k_1}=( a^{k_2} b^{(l-1)})b,$$ 
if $l> 1$, or the equality
$$  a^{k_1}= a^{k_2}b,$$ 
 if $l=1.$  

Finally, the equality $a^n=bu$ implies 
 the equality $ a^{n+1}=bua$, and hence \eqref{suk} holds for $x=a^{n}$ and $y=ua$. The proof in the ``left'' case is similar. 
\qed 

\vskip 0.2cm

The above definitions imply that a semigroup $S$ of rational functions is power joined if 
for all $A,B\in S$ there exist $k,l\geq 1$ such that 
 \be \l{eba} A^{\circ k}=B^{\circ l},\ee
or in other words if any two elements of $S$ share an iterate.  On the other hand, $S$  is  
power twisted if for all $A,B\in S$ there exist $k_1,k_2,l\geq 1$ such that
\be \l{piz0} A^{\circ k_1}=A^{\circ k_2}\circ B^{\circ l}.\ee
 Finally, $S$ is left (resp. right) Archimedean if for all $A,B\in S$ the function $B$ is a compositional right (resp. left) factor of some iterate of $A$. 
 Notice  that 
if a semigroup $S$ of rational functions is power joined or power twisted and $\underline{S}$ is not empty, then $S=\underline{S}$, since for a rational function $A$ such that $\deg A=1$ any of equalities \eqref{eba}, \eqref{piz0} implies that $\deg B=1$.

For a semigroup of rational functions $S$ the condition that $S$ is power twisted can be 
replaced by an apparently stronger condition, which naturally arises in the study of rational functions sharing a measure of maximal entropy. Namely, the following statement holds. 

\bl \l{ptw} Let $S$ be a semigroup of rational functions of degree at least two. Then $S$ 
is power twisted if and only if 
for all $A,B\in S$ there exist $k,l\geq 1$ such that the equalities 
\be \l{piz01} A^{\circ 2k}=A^{\circ k}\circ B^{\circ l}, \ \ \ \ \  B^{\circ 2l}=B^{\circ l}\circ A^{\circ l}
\ee hold. 
\el
\pr The ``if'' part is obvious. To prove the ``only if'' part, let us show first that equality \eqref{piz0}  implies that there exists $k,l\geq 1$ such that the first equality in \eqref{piz01} holds. Comparing degrees in \eqref{piz0}, we see that $k_1>k_2$. Therefore \eqref{piz0} can be rewritten in the form  
\be \l{skun} A^{\circ s}\circ A^{\circ k}=A^{\circ s}\circ B^{\circ l}\ee for some $s,k,l\geq 1$.  
Clearly, \eqref{skun}   
implies that 
\be \l{puk} A^{\circ (s+t)}\circ A^{\circ k}=A^{\circ (s+t)}\circ B^{\circ l}\ee for every $t\geq 0.$  
If $k-s\geq 0$, then setting $t=k-s$ in \eqref{puk}, we obtain the needed equality.  The general case reduces to this one, since \eqref{skun} implies that for every $r\geq 1$ the 
equality 
\be \l{pak} A^{\circ s}\circ A^{\circ kr}=A^{\circ s}\circ B^{\circ lr}\ee holds. Thus, for $r$ big enough,  
\eqref{skun} holds for $k'=kr$ and $l'=lr$ with $k'-s\geq 0.$ 

The above shows that for all $A,B\in S$ there exist 
 $k_1,l_1\geq 1$ and $k_2,l_2\geq 1$ such that the equalities  
\be \l{pzz01} A^{\circ 2k_1}=A^{\circ k_1}\circ B^{\circ l_1}, \ \ \ \ \ \ B^{\circ 2l_2}=B^{\circ l_2}\circ A^{\circ k_2}\ee hold. Moreover, since \eqref{pzz01} implies that  
$$(\deg A)^{k_1}=(\deg B)^{l_1}, \ \ \ \ \  (\deg A)^{k_2}=(\deg B)^{l_2},$$ the equality $l_1k_2=l_2k_1$ holds.  

Since equalities \eqref{pzz01} imply that for all $s,r\geq 1$ 
the equalities  
\be \l{pzz02} A^{\circ 2k_1s}=A^{\circ k_1s}\circ B^{\circ l_1s}, \ \ \ \ \ \ B^{\circ 2l_2r}=B^{\circ l_2r}\circ A^{\circ k_2r}\ee hold, setting 
$s=k_2$ and $r=k_1,$ we see that the equalities \eqref{piz01}
 hold for $k=k_1k_2$ and $l=l_1k_2=l_2k_1.$ \qed 

\bl \l{svin2} Let $S$ be a semigroup of rational functions such that $\underline{S}$ is finite. If $\overline{S}$  is power joined, then $S$ is  left and right reversible. On the other hand, if $\overline{S}$  is power twisted, then $S$ is  right reversible.
\el
\pr  If $\overline{S}$  is power joined, then $\overline{S}$ is left and right reversible by Lemma \ref{svin1}. Thus, to prove the first part of the lemma we only must construct solutions $X,Y\in S$ of the equations 
\be \l{bive2} A\circ X=B\circ Y \ \ \ \ {\rm and} \ \ \ \ X\circ A=Y\circ B\ee in case if at least one of the functions $A,B\in S$ is of degree one.

Assume, say, that $\deg A=1$. Then it follows from the finiteness of $\underline{S}$  that 
\be \l{ife} A^{\circ k}=z\ee
for some $k\geq 1$, implying in particular that  the function $z$ belongs  to $S$.
 Therefore, 
equalities \eqref{bive2} hold for $X,Y\in S$ given by the formulas 
$$X= A^{\circ {(k-1)}}\circ B, \ \ \ \ Y=z,\ \ \ \ {\rm and} \ \ \ \ X= B\circ A^{\circ {(k-1)}}, \ \ \ \ Y=z,$$ correspondingly (in case $k=1$ we set $A^{\circ {(k-1)}}=z$).

Similarly, if $\overline{S}$  is power twisted, then $\overline{S}$ is right reversible by Lemma \ref{svin1}, and, assuming that \eqref{ife} holds, we see that   
 the second equality in \eqref{bive2} is satisfied for 
$$X= B\circ A^{\circ {(k-1)}}, \ \ \ \ Y=z. \eqno{\Box} $$

\subsection{Semigroups $C_{\infty}(P)$}


We will call a rational function  {\it special} if it is 
either a Latt\`es map or is conjugate to  $z^{\pm n},$ $n\geq 2$, or $\pm T_n,$ $n\geq 2$, where $T_n$ is the $n$th Chebyshev polynomial. Since a polynomial cannot be a Latt\`es map, this definition is consistent with the definition of special polynomials given in the introduction.  

For a rational function of degree at least two, we denote by $C(P)$  
the collection of rational functions, including rational functions of degree one, commuting with $P$. It is clear that $C(P)$ is  a semigroup. For  
the subsemigroup $\underline{C(P)}$ of $C(P)$ we will use the standard notation  $\Aut(P).$ 
It is easy to see that $\Aut(P)$ is a group. Moreover, since elements of $\Aut(P)$ permute 
fixed points of $P^{\circ k},$ $k\geq 1$, and any M\"obius
transformation is defined by its values at any three points, the group $\Aut(P)$ is finite. 
In particular, for every $A\in \Aut(P)$ equality \eqref{ife}   
holds for some $k\geq 1$. 

The following fact is proved easily by a direct calculation (see \cite{rev}, Lemma 2.1).

\bl \l{per} If $A$ and $U$ are rational functions such that 
$A\circ U\in C(P)$ and $U\in C(P),$ then  $A\in C(P)$. \qed 
\el

Let us define the sets $C_{\infty}(P)$  and $\Aut_{\infty}(P)$ by the formulas 
\be \l{defin} C_{\infty}(P)= \bigcup_{i=1}^{\infty}C(P^{\circ k}), \ \ \ \ \Aut_{\infty}(P)=\bigcup_{k=1}^{\infty} \Aut(P^{\circ k}).\ee
Since obviously  
\be \l{lcm} C(P^{\circ k}),\  C(P^{\circ l})\subseteq C(P^{\circ \LCM(k,l)})\ee 
and 
  \be \l{xre} \Aut(P^{\circ k}),\   \Aut(P^{\circ l})\subseteq  \Aut(P^{\circ \LCM(k,l)}),\ee 
the set $C_{\infty}(P)$ is a semigroup, and the set  $\Aut_{\infty}(P)$ is a group.

\bl \l{zaqw0} Let $P$ be a polynomial of degree at least two, and $S$  a semigroup of polynomials such that the semigroup $\overline{S}$ is non-empty. Then $S$
 is  contained in  $C(P)$ if and only if $\overline S$ is  contained in  $\overline{C(P)}$. 
Similarly,  $S$
 is  contained in  $C_{\infty}(P)$ if and only if $\overline S$ is  contained in  $\overline{C_{\infty}(P)}$. 
\el 
\pr The ``only if'' parts of the lemma are  clear. To prove the ``if'' parts, we observe that if $Q$ is any fixed element of $\overline{S}$, then for every $\alpha\in\underline{S}$ the function  $\alpha\circ Q$ belongs to $\overline{S}$. Therefore, if  
$\overline S$ is  contained in  $\overline{C(P)}$, then both $Q$ and $\alpha\circ Q$ belong to $C(P)$, implying  by Lemma \ref{per} that $\alpha$ belongs to $\Aut(P)$. Thus, $\underline{S}\subseteq C(P).$ 
 Similarly, if $\overline S$ is  contained in  $\overline{C_{\infty}(P)}$, then   $\alpha\circ Q\in C(P^{\circ k})$ for some $k\geq 1$, implying that  $\alpha\in\Aut(P^{\circ k})$. \qed 

\vskip 0.2cm
We recall that by the Ritt theorem (see \cite{r} and also \cite{e2}, \cite{rev}) if rational functions $A$ and $B$ of degree at least two commute, then either they both are special or they have an iterate in common.
This result implies the following characterizations of semigroups $\overline{ C_{\infty}(P)}$.

\bl \l{svi} For every non-special rational function $P$ of degree at least two the semigroup $\overline{C_{\infty}(P)}$ coincides with the set of rational functions sharing an iterate with $P.$ 
\el
\pr If $A$ commutes with some iterate of $P$, then the Ritt theorem implies that $A$ and $P$ share an iterate.  On the other hand, if there exist $k,l\in \N$ such that $A^{\circ k}=P^{\circ l},$ then $A$ obviously 
commutes with $P^{\circ l}.$ \qed 

\vskip 0.2cm

In turn, subsemigroups of $C_{\infty}(P)$ can be characterized as follows.

\bt \l{3}  Let $S$ be a power joined  
semigroup of rational functions of degree at least two.  Then $S$ is a subsemigroup
of the semigroup $\overline{C_{\infty}(P)}$ for every $P\in S$. In the other direction,  every subsemigroup
of the semigroup $\overline{C_{\infty}(P)}$, where $P$ is a non-special rational function of degree at least two, is power joined.
\et
\pr If $S$ is power joined and $P\in S$, then for an arbitrary element  $A\in S$ the equality $A^{\circ k}=P^{\circ l}$ holds for some $k,l\in \N$,  implying that $A$ commutes with $P^{\circ l}.$ Therefore, $$S\subseteq \bigcup_{i=1}^{\infty}C(P^{\circ k})=C_{\infty}(P).$$

In the other direction, if  $S\subseteq C_{\infty}(P)$, 
 then \eqref{lcm} implies that for all $A,B\in S$ there exist $l\in \N$ such that both $A$ and $B$ commute with  $P^{\circ l}$. It follows now from the Ritt theorem   that there exist $k_1,k_2,r_1,r_2\in \N$ such that the equalities
$$A^{\circ k_1}=P^{\circ lr_1}, \ \ \ \  \ B^{\circ k_2}=P^{\circ lr_2}$$
hold, implying that 
 $$A^{\circ k_1r_2}=B^{\circ k_2r_1}.\eqno{\Box}$$

\bt \l{1} Let  $P$ be a  non-special rational function of degree at least two.
Then every subsemigroup $S$ 
of the semigroup $C_{\infty}(P)$ is cancellative and left and right  reversible. 
\et 
\pr  Suppose that  
 \be\l{na} F\circ X=F \circ Y\ee for some $F,X,Y\in C_{\infty}(P)$. Clearly, if $\deg F=1,$ then $X=Y,$ so we can assume that $\deg F>1.$ 
Let $k,l\in \N$ be numbers such that  $F^{\circ k}=P^{\circ l}$, and  $s\in \N$ a number such that both $X,Y$ commute with $P^{\circ s}.$ 

Obviously, both $X,Y$ commute  with $P^{\circ ls}=F^{\circ ks}.$ 
Since equality  
\eqref{na} implies the equality 
$$F^{\circ ks}\circ X=F^{\circ ks} \circ Y,$$ this yields that 
$$X\circ F^{\circ ks}=Y\circ F^{\circ ks},$$ implying that $X=Y.$  
Therefore, the semigroup $C_{\infty}(P)$ is cancellative, implying that every its subsemigroup is also cancellative. 

Further, since $\overline{C_{\infty}(P)}$ is power joined by Theorem \ref{3} and $\Aut(P)=\underline{C_{\infty}(P)}$ is finite, for every  subsemigroup $S$ of $C_{\infty}(P)$ the semigroup $$\overline{S}=S\cap \overline{C_{\infty}(P)}$$ is power joined, and the semigroup  $$\underline{S}=S\cap \underline{C_{\infty}(P)}$$ is finite. Thus,  
the left and the right reversibility of $S$ follow from Lemma \ref{svin2}. \qed

\subsection{Amenability of semigroups $C_{\infty}(P)$}

Let $P$ be a non-special rational function of degree at least two. Following \cite{rev}, we define an equivalence relation 
$\underset{P}{\sim}$ on the semigroup $C(P)$, setting  
$Q_1\underset{P}{\sim} Q_2$ if 
\be \l{zek} Q_1\circ P^{\circ l_1}=Q_2\circ P^{\circ l_2}\ee  
for some $l_1\geq 0,$ $l_2\geq 0$. 

The following lemma is an easy corollary of the right cancellativity of semigroups of rational functions 
(see \cite{rev}, Lemma 3.1).

\bl \l{eg}  Let $\bf{A}$ be an equivalence class of $\underset{P}{\sim}$.  
For any $n\geq 1$ the class $\bf{A}$   contains at most one rational function of degree $n$. Furthermore, if $A_0\in \bf{A}$ is a function of minimum possible degree, then 
any $A\in \bf{A}$ has the form $A=A_0\circ P^{\circ l},$ $l\geq 0.$  \qed 
\el  

The following result was proved in \cite{rev}.

\bt \l{gr} Let $P$ be a non-special rational function of degree at least two. Then the relation $\underset{P}{\sim}$ is a congruence on the semigroup $C(P)$, and the quotient semigroup is a finite  group. 
\qed 
\et

It was shown in the paper \cite{peter} that 
every power joined subsemigroup of rational functions is amenable.
Below we reprove this result in a slightly more general form. As in the paper \cite{peter}, our proof relies on Theorem \ref{gr}. However, our reduction to Theorem \ref{gr} is different and uses  
 the F{\o}lner criterion.

\bt \l{2}  Let  $P$ be a  non-special rational function of degree at least two.
Then every subsemigroup $S$ 
of the semigroup $C_{\infty}(P)$ is   amenable. 
\et 
\pr
By Theorem \ref{1}, the semigroup $C_{\infty}(P)$ is cancellative. Furthermore,  $C_{\infty}(P)$ cannot contain a free subsemigroup on two generators. Indeed, if $A,B\in S$ are of degree greater than one, then $A$ and $B$ have a common iterate and hence $\langle A,B\rangle$ is not free. On the other hand, if say $A$ is of degree one, then $\langle A,B\rangle$ is not free since \eqref{ife} implies that $A^{\circ (k+1)}=A$. 
Therefore, by Theorem \ref{fr}, to prove the theorem we only must show that  $C_{\infty}(P)$ is amenable. Moreover, it follows from \eqref{lcm} by 
Lemma \ref{da} that it is enough to prove the  amenability of the semigroups $C(P^{\circ k})$, $k\geq 1$. Finally, since iterates of a non-special rational function $P$ are non-special 
(see \cite{fin}, Lemma 2.12),  it is enough to  prove only the  amenability of $C(P)$ for an arbitrary non-special rational function $P$.

By Lemma \ref{eg} and Theorem \ref{gr}, there exist  $X_1,X_2,\dots , X_n\in C(P)$ such that $$C(P)=\bigsqcup_{i=1}^nM_i,$$ where 
$$M_i=\{X\in C(P)\, \vert \, X=X_i\circ P^j,\ j\geq 0\}.$$
For $N\geq 0$ and $i,$ $1\leq i \leq n$, we set
$$M_{i,N}=\{X\in C(P)\, \vert \, X=X_i\circ P^j,\, 0\leq j\leq N\}$$
and $$F_N=\bigsqcup_{i=1}^nM_{i,N}.$$ Let us show that for 
every finite subset $H$  of $C(P)$ and every $\varepsilon >0$ the set $F_N$ with $N$ big enough satisfies the condition 
\be \l{t} \vert F_N\setminus X\circ F_N\vert \leq \varepsilon \vert F_N\vert\ee
for all $X\in H.$ 

By Theorem \ref{gr}, for every $j,i,$ $1\leq j,i \leq n$ there exist $m(j,i)\in \N$ and $k(j,i),$ $1\leq k(j,i) \leq n,$ such that 
\be \l{gop1} X_j\circ X_i=X_{k(j,i)}\circ P^{\circ m(i,j)}.\ee Moreover, for fixed $j$ the map $i\rightarrow k(j,i)$ is a
bijection of the set $\{1,2,\dots,n\}.$
Set  $$L_1=\max_{1\leq j,i\leq n}m(j,i).$$
Since $H$ is a subset of $C(P)$, every element $X$ of $H$ can be represented in the form  
\be \l{gop2} X=X_j\circ P^{\circ l}\ee for some $j,$ $1\leq j \leq n,$ and $l\geq 0,$ and we define $L_2$ as the maximum number $l$ in such a representation (since $H$ is finite, such a number exists). Clearly, $$\vert F_N\vert =(N+1)n,$$ and it follows from \eqref{gop1} and \eqref{gop2} that for every $X\in H$ the inequality 
$$\vert F_N\setminus X\circ F_N\vert \leq n(L_1+L_2)$$ holds. 
Therefore,  \eqref{t} holds for $N$ big enough and hence $C(P)$ is left amenable by Theorem \ref{fo}. Since the set $M_i,$ $1\leq j \leq n,$ coincides with the set
$$M_i'=\{X\in S\, \vert \, X=P^j\circ X_i,\ j\geq 0\},$$
a symmetric argument shows that $C(P)$ is right amenable.
\qed

\section{Semigroups $E(P)$}
Let us recall that for a  rational function $P$ of degree at least two we denote by $\mu_P$ the measure of 
 maximal entropy for $P$, and by $E(P)$ the set of rational functions $Q$ of degree at least two such that $\mu_Q=\mu_P$, completed by  $\mu_P$-invariant M\"obius transformations.

\bl \l{kosh1} Let $P$ be a rational function of degree at least two. Then the set $E(P)$ is a semigroup. 
\el 
\pr 
Let $A$ and $B$ be elements of $E(P)$ of degree $n$ and $m$ correspondingly.  Assume first that $n,m\geq 2.$ We recall that the measure $\mu_P$ is characterized by the balancedness property that $$\mu_P(P(S))=\mu_P(S)\deg P\,$$ 
for any Borel set $S$ on which $P$ is injective (\cite{flm}). Therefore, we only must show that if $\mu_P$ is 
the balanced measure for $A$ and $B$, then $\mu_P$ is 
the balanced measure for $A\circ B$. 
Let $S$ be a Borel set on which $A\circ B$ 
is injective. Then $B$ is injective on $S$ and $A$ is injective on $B(S)$, implying that $$\mu_P\big((A\circ B)(S)\big)=\mu_P\big(A(B(S)\big)
=n\mu_P\big(B(S)\big)=nm\mu_P(S).$$
Thus, $\mu_P$ is the balanced measure for $A\circ B$.

Further, if $A\in E(P)$ is a function of degree $n\geq 2$, and  $\sigma$ is a $\mu_P$-invariant  M\"obius transformation, then for any  Borel set $S$  on which $A\circ \sigma$ is injective we have $$\mu_P\big((A\circ \sigma)(S)\big)=\mu_P\big(A(\sigma(S)\big)
=n\mu_P\big(\sigma(S)\big)=n\mu_P(S).$$ 
Similarly, for any  Borel set $S$  on which $\sigma\circ A$ is injective we have
 $$\mu_P\big((\sigma\circ A)(S)\big)=\mu_P\big(\sigma(A(S)\big)
=\mu_P\big(A(S)\big)=n\mu_P(S).$$ 
 Thus,   
$\mu_P$ is the balanced measure for $A\circ \sigma$ and $\sigma \circ A.$ 

Finally, it is clear that if  $\sigma_1$ and $\sigma_2$ are $\mu_P$-invariant  M\"obius transformation, then $\sigma_1\circ \sigma_2$ is also such a transformation. \qed 

\vskip 0.2cm

Algebraic conditions for non-special rational functions $A$ and $B$ to share a measure of maximal entropy were obtained in the papers \cite{lev}, \cite{lp}, and can be formulated as 
follows (see \cite{ye}).

\bt \l{ye} Let $A$ and $B$ be non-special rational functions of degree at least two. Then $\mu_A=\mu_B$ if and only if there exist $k,l\geq 1$ such that the equalities 
\be \l{sisis} A^{\circ 2k}=A^{\circ k}\circ B^{\circ l},\ \ \ \ \ \ \ \ \ B^{\circ 2l}=B^{\circ l}\circ A^{\circ l},\ee
hold. \qed 
\et 

Notice that either of 
equalities in \eqref{sisis} is sufficient for the equality $\mu_A=\mu_B$, regardless whether $A$ and $B$ are special or not (\cite{lev}).  Notice also that in a sense describing solutions of the system \eqref{sisis} reduces to describing rational functions which are not tame (see \cite{ye}, \cite{entr}).

Rational functions sharing an iterate share a measure of maximal entropy, and the system \eqref{sisis} can be regarded as a generalization of the condition that $A$ and $B$ share an iterate. Correspondingly,  the following statement takes the place of Theorem \ref{3}.

\bt \l{kosh2} Let $S$ be a power twisted 
semigroup of rational functions of degree at least two.  Then $S$ is a subsemigroup
of the semigroup $\overline{E(P)}$ for every $P\in S$. In the other direction,  every subsemigroup
of the semigroup $\overline{E(P)}$, where $P$ is a non-special rational function of degree at least two, is power twisted.
\et
\pr If $S$ is power twisted, and $P\in S$, then by Lemma \ref{ptw} for every $A\in S$ there exist $k$ and $l$ such that $A^{\circ 2k}=A^{\circ k}\circ P^{\circ l}$, implying that $A\in \overline{E(P)}$.   

In the other direction, it is well known that if $P$ is non-special, then all rational functions sharing a measure of maximal entropy with $P$ also are non-special. Moreover, if $A,B$ are such functions, then  by Theorem \ref{ye} equalities \eqref{sisis} hold, implying that every  subsemigroup of $\overline{E(P)}$ is power twisted. \qed  

Finally, the following lemma is the analogue of Lemma \ref{zaqw0}.  

\bl \l{zaqw} Let $P$ be a polynomial of degree at least two, and $S$  a semigroup of polynomials such that the semigroup $\overline{S}$ is non-empty. Then $S$
 is  contained in  $E(P)$ if and only if $\overline S$ is  contained in  $\overline{E(P)}$. 
\el 
\pr The ``only if'' part is clear. To prove the ``if'' part, 
we observe that if $Q$ is any fixed element of $\overline{S}$, then for every $\alpha\in\underline{S}$ the function $Q\circ \alpha\in \overline{S}$  belongs to $\overline{S}$. 
 Thus, by the invariance of $\mu_P$, for any Borel set $S$ we have:
 $$\mu_P\big((Q\circ \alpha)^{-1}(S)\big)=\mu_P(S).$$ On the other hand,  
$$\mu_P\big((Q\circ \alpha)^{-1}(S)\big)=\mu_P\big(Q^{-1}(\alpha^{-1}(S)\big)=
\mu_P\big(\alpha^{-1}(S)\big).$$
Therefore, 
$$\mu_P\big(\alpha^{-1}(S)\big)=\mu_P(S),$$ implying that $\alpha \in \underline{E(P)}$. \qed

\section{Semidirect products}\l{ameb}


Let us recall that for a  rational function $R$ of degree at least two, the group $ G(R)$ is defined  as the group of M\"obius transformations $\sigma$ such that \be \l{eblys} R\circ \sigma=\nu \circ R\ee for some  M\"obius transformations $\nu$. 
 It is easy to see that  $ G(R)$ is indeed a group 
and that the map \be \l{homic}\gamma_R:\sigma \rightarrow \nu_{\sigma}\ee is a homomorphism from   $ G(R)$ to the group $\rm{Aut}(\C\P^1)$. Notice  that the group ${\rm Aut}(R)$ is a subgroup of  $ G(R)$. We say that a rational function $R$ of degree $n\geq 2$ is {\it a quasi-power} if there exist $\alpha,\beta\in \Aut(\C\R^1)$ such that $$R= \alpha\circ z^n\circ \beta.$$

The following statement was proved in \cite{fin} (see also \cite{sym} for more results about  $G(R)$ and related groups). 

\bt \l{prim} Let $R$ be a rational function of degree at least two that is not a quasi-power. Then the group 
 $ G(R)$ is finite. \qed
\et 


Assume that $\Gamma$ a subgroup of $G(R)$ such that $\gamma_R(\Gamma)\subseteq \Gamma$. Then \eqref{homic} is an endomorphism of $G(R)$. Furthermore, $\gamma_R$ defines in an obvious way a homomorphism \be \l{moh} \rho_R: \langle R\rangle \rightarrow \rm{End}(\Gamma).\ee  We denote by $S_{\Gamma,R}$ the semigroup generated by $\Gamma$ and $R$. It is clear that a rational function $A$ belongs to $S_{\Gamma,R}$ if and only if 
 \be \l{form0} A=\delta \circ R^{\circ s}\ee for some $s\geq 0$ and $\delta\in \Gamma.$ Moreover, in the notation of Section \ref{xcv}, we have: 
$$S_{\Gamma,R}=\Gamma\underset{\rho_{R}}{\times} \langle R\rangle.$$

\bl \l{fe}
Let $R$ be a  rational function of degree $n\geq 2$  not conjugate to $z^{\pm n}$,  and $\Gamma$ a subgroup of $G(R)$ such that 
$\gamma_R(\Gamma)\subseteq \Gamma$. Then $\Gamma$ is finite.  
\el
\pr By Theorem \ref{prim}, $G(R)$ and hence $\Gamma$  is finite, unless $R$ is a quasi-power.  On the other hand, since the group $G(z^n)$ consists of the M\"obius transformations $cz^{\pm 1},$ $c\in \C\setminus\{0\}$ (see \cite{fin}, Lemma 4.1), it is easy to see that the condition $\gamma_R(\Gamma)\subseteq \Gamma$ holds for a quasi-power $R$ only if $R$  is conjugate to $z^{\pm n}.$  \qed 

\vskip 0.2cm

The following two results describe amenability properties of subsemigroups of $S_{\Gamma,R}$ according to whether the condition  $\gamma_R(\Gamma)= \Gamma$ or the condition $\gamma_R(\Gamma)\subseteq \Gamma$ is satisfied. In particular, they imply Theorem \ref{pi} from the introduction.

\bt  \l{moh2} Let $R$ be a  rational function of degree $n\geq 2$  not conjugate to $z^{\pm n}$,  and $\Gamma$ a subgroup of $G(R)$ such that 
$\gamma_R(\Gamma)\subseteq \Gamma$. Then $S_{\Gamma,R}$ is left amenable if and only if
  $\gamma_R(\Gamma)= \Gamma$. Moreover, if $\gamma_R(\Gamma)= \Gamma$, then $S_{\Gamma,R}\subseteq C_{\infty}(R)$ and every subsemigroup of   $S_{\Gamma,R}$ is amenable. 
\et 
\pr 
Since all elements of $S_{\Gamma,R}$ have the form \eqref{form0}, if $\sigma_0\in \Gamma$ but $\sigma_0\not\in {\rm Im}\, \gamma_R$, the equation 
$$R\circ X=(\sigma_0\circ R)\circ Y$$ has no 
 solutions $X,Y$ in $S_{\Gamma,R}$.  Therefore, whenever $\gamma_R(\Gamma)$ is a proper subset of $\Gamma$, the semigroup  
$S_{\Gamma,R}$ is not left reversible and hence is not left amenable. 

On the other hand, since the group $\Gamma$ is finite by Lemma \ref{fe}, if $\gamma_R(\Gamma)= \Gamma$, then the restriction $\gamma_R:\Gamma \rightarrow \Gamma$ is an automorphism. Moreover, since the automorphism group of a finite group is finite,  there exists $l\geq 1$ such that the  iterate $\gamma_R^{\circ l}$ is the identical automorphism.
Therefore, since 
$$R^{\circ l}\circ \sigma=\gamma_R^{\circ l}(\sigma)\circ R^{\circ l}, \ \ \ \ \sigma\in \Gamma,$$  the group $\Gamma$ is a subgroup of $\Aut(R^{\circ l})$.  
In turn, this implies that 
for every element $A\in \overline{S_{\Gamma,R}}$ the iterate $A^{\circ l}$ commute with $R^{\circ l}$, since 
$$A^{\circ l}=\sigma \circ R^{\circ sl}$$ for some $\sigma\in \Gamma$ and $s\geq 1$ by \eqref{form0}.  Hence, by the Ritt theorem, $A$ and $R$ share an iterate. Therefore,
$\overline{S_{\Gamma,R}}\subseteq \overline{C_{\infty}(R)}$  by Lemma \ref{svi}, 
implying that $S_{\Gamma,R}\subseteq C_{\infty}(R)$ by Lemma \ref{zaqw0}. 
In particular, if  $S'$ is a subsemigroup of $S$, then $S'$ is a subsemigroup of $C_{\infty}(R)$, implying that $S'$  is amenable by Theorem \ref{2}.  \qed

\bt \l{moh1} Let $R$ be a  rational function of degree $n\geq 2$  not conjugate to $z^{\pm n}$, and $\Gamma$ a subgroup of $G(R)$ such that 
$\gamma_R(\Gamma)\subseteq \Gamma$. Then $S_{\Gamma,R}\subseteq E(R)$ and every subsemigroup of $S_{\Gamma,R}$ is  right amenable. 
\et 
\pr By Lemma \ref{zaqw} and Theorem \ref{kosh2}, to prove that $S_{\Gamma,R}\subseteq E(R)$ it is enough to show that 
 $\overline{S_{\Gamma,R}}$ is power twisted, that is, that for all $A,B\in \overline{S_{\Gamma,R}}$ there exist $k,l\geq 1$ such that \eqref{piz0} holds. It follows from the representation \eqref{form0} that considering instead of $A$ and $B$ some of their  iterates without loss of generality we may assume that $\deg A=\deg B$ and \be \l{bru} A=\sigma\circ B\ee for some $\sigma \in \Gamma$. Furthermore, it follows from \eqref{form0} and \eqref{bru}   that for every $k\geq 1$ there exists $\sigma_k\in \Gamma$ such that $$A^{\circ k}=\sigma_k\circ B^{\circ k}.$$  Therefore, since $\Gamma$ is 
finite, there exist  $k_1,k_2\geq 1$ such that $k_1>k_2$ and 
$$A^{\circ k_1}=\delta\circ B^{\circ k_1}, \ \ \ \  A^{\circ k_2}=\delta\circ B^{\circ k_2}$$ for the same  $\delta \in \Gamma,$ implying that   \eqref{piz0} holds for $l=k_1-k_2.$

Set 
$$\Gamma_k=\gamma_R^{\circ k}(\Gamma).$$ Since $\Gamma$ is finite, it follows from $$\Gamma\supseteq \Gamma_1 \supseteq \Gamma_2 \supseteq \dots  $$ that there exists $k_0$ such that $ \Gamma_k = \Gamma_{k_0}$ for all $k\geq k_0$.  We set $\widehat \Gamma=\Gamma_{k_0}$ and $$\Gamma_0={\rm Ker}\,\gamma_R^{\circ k_0},$$ so that 
$\widehat \Gamma=\Gamma/\Gamma_0.$ 
Since $\gamma_R:\widehat \Gamma \rightarrow \widehat \Gamma$ is an isomorphism, the above definitions imply that for $\alpha_1,\alpha_2\in \Gamma$  the equality 
\be \l{eas} \gamma_R^{\circ k}(\alpha_1)= \gamma_R^{\circ k}(\alpha_2)\ee holds for some $k\geq k_0$ if and only if elements $\alpha_1$ and $\alpha_2$ belong to the same coset of $\Gamma_0$ in $\Gamma$. 

Since $\underline{S_{\Gamma,R}}=\Gamma$ is finite and  $\overline{S_{\Gamma,R}}$ is power twisted, $S_{\Gamma,R}$ is right reversible by Lemma \ref{svin2}. Thus, by the first part of Theorem \ref{lastik}, equivalence classes on $S_{\Gamma,R}$ corresponding to equivalence relation \eqref{listik} form a semigroup $S_{\Gamma,R}/\sim $. Let us show that $S_{\Gamma,R}/\sim $ is isomorphic to the semigroup $S_{\widehat\Gamma,R}.$ For this purpose, it is enough to prove the following statement: 
for $\alpha_1,\alpha_2\in \Gamma$ and $s_1,s_2\geq 0$
the equality 
\be \l{mir} (\beta \circ R^{\circ s})\circ (\alpha_1\circ R^{\circ s_1})=(\beta \circ R^{\circ s})\circ (\alpha_2\circ R^{\circ s_2})\ee  holds for some $s\geq 0$ and  $\beta\in \Gamma$ if and only if $s_1=s_2$ and $\alpha_1,\alpha_2$ belong to the same coset of $\Gamma_0$ in $\Gamma$. To prove the ``if'' part, we observe that if $\alpha_2=\delta\circ \alpha_1,$ where $\delta \in \Gamma_0$, then \eqref{mir} holds 
for $s=k_0$ and any $\beta.$ On the other hand, if equality \eqref{mir} holds, then obviously $s_1=s_2$ and 
\be \l{miru} R^{\circ s}\circ (\alpha_1\circ R^{\circ s_1})= R^{\circ s}\circ (\alpha_2\circ R^{\circ s_2}), 
\ee implying that for every $l\geq 0$ the equality 
$$R^{\circ (s+l)}\circ (\alpha_1\circ R^{\circ s_1})= R^{\circ (s+l)}\circ (\alpha_2\circ R^{\circ s_2})$$ holds. Thus, without loss of generality we may assume that $s\geq k_0$ in \eqref{miru}, implying that $\alpha_1$ and $\alpha_2$ belong to the same coset of $\Gamma_0$ in $\Gamma$.

Since $S_{\Gamma,R}/\sim $ is isomorphic to the semigroup $S_{\widehat\Gamma,R},$ it follows from the second part of Theorem \ref{lastik} that to prove that any subsemigroup of $S_{\Gamma,R}$ is right amenable  it is enough to prove that 
any subsemigroup of $S_{\widehat\Gamma,R}$ is  right amenable. In turn, the last statement follows from Theorem \ref{moh2}, which implies that $S_{\widehat\Gamma,R}$ is  amenable.  \qed

\vskip 0.2cm

Notice that the proofs of Theorem \ref{moh2} and Theorem \ref{moh1} remain true for $R$ conjugate to $z^{\pm n}$  if to require the finiteness of $\Gamma$. For example, the semigroup generated by the polynomial $z^2$ and the M\"obius transformation $z\rightarrow - z$ is the simplest example of a  right amenable semigroup of rational functions that is not left amenable. On the other hand, the semigroup generated by $z^3$ and the M\"obius transformation $z\rightarrow e^{\frac{2\pi i}{5}}z$, say,  is amenable.

\section{Tame semigroups of rational functions}
\subsection{Tame rational functions}
We recall that a rational function $A$ of degree at least two is called {\it tame} if the algebraic curve $$A(x)-A(y)=0$$ has no factors of genus zero or one distinct from the diagonal. 
By the Picard theorem, this condition  is equivalent to the condition that the equality   
\be \l{urav} A\circ f=A\circ g,\ee where  $f$ and $g$  are functions  meromorphic on $\C$, implies that $f\equiv g.$ 
 Notice that any rational function of degree two  is not tame since the curve 
$$\frac{A(x)-A(y)}{x-y}=0$$ has degree one, implying that its genus is zero. Thus, a tame rational function has degree at least three. 
Notice that a {\it general} rational function of degree at least four is tame.  Specifically, a rational function of degree at least four is tame whenever it has only simple critical values (\cite{pg}). 

We say that a semigroup of rational functions $S$ is {\it tame}, if it contains tame rational functions only. Clearly, the  tameness condition can be regarded as 
a  strengthening of the cancellativity condition.
\bl \l{ltam} Tame rational functions form a cancellative subsemigroup of $\C(z).$
\el 
\pr Let us assume that $A$, $B$ are tame rational functions and  $f$, $g$ are   meromorphic on $\C$ functions such that the equality 
$$ (B\circ A)\circ f=(B\circ A) \circ g$$ holds.  Since  $B$ is tame and $A\circ g$ and $A\circ f$  are meromorphic on $\C$, this equality implies equality \eqref{urav}. In turn,  equality  \eqref{urav} implies  that $f\equiv g$, since $A$ is tame. Thus, tame rational functions form a subsemigroup of $\C(z)$, and it is clear that this subsemigroup is cancellative.  
\qed 

\vskip 0.2cm

Our approach to the amenability of tame semigroups of rational functions is based on the three results about tame rational functions  from the paper \cite{tame} given below. 

Let $P^{\circ d}= U\circ V$ be a decomposition of an iterate $P^{\circ d}$ of a rational function $P$ into a composition of rational functions $U$ and $V$. 
We say that this decomposition is {\it induced} by a decomposition $P^{\circ d'}= U'\circ V'$, where $d'<d,$ if there exist 
$k_1,k_2\geq 0$ such that 
$$U=P^{\circ k_1}\circ U', \ \ \ \ \  V=V'\circ P^{\circ k_2}.$$  
The first statement we need is following (\cite{tame}). 

\bt \l{1+} Let $P$ be a tame rational function of degree $n$. Then there exists an integer $N$, depending on  $n$ only, such that  any decomposition of  $P^{\circ d}$ with $d\geq N$ is induced by a decomposition of $P^{\circ N}$. \qed 
\et

We recall that functional decompositions  $R=U\circ V$ of a rational function $R$ into compositions of rational functions 
$U$ and $V$,  considered up to the equivalence 
\be \l{equi} U\rightarrow U\circ \alpha, \ \ \ \ V\rightarrow\alpha^{-1} \circ V,\ \ \ \ \alpha\in Aut(\C\P^1),\ee are in a 
one-to-one correspondence with imprimitivity systems of the monodromy group of $R$. In particular, the number of such classes is finite. Consequently, Theorem \ref{1+} implies that 
for every tame rational function $P$ there exist finitely many rational functions $F_1,F_2,\dots, F_s$ such that a rational function $F$ is a compositional right factor of an iterate of $P$ if and only if  $F$ has the form \be \l{egik} F=\alpha\circ F_i\circ P^{\circ l}, \ \ \ l\geq 0,  \ \ \ 1\leq i \leq s,\ \ \ \alpha\in Aut(\C\P^1).\ee

It is easy to see that if rational functions $A$ and $B$ have a common iterate, then 
each iterate of $B$ is a compositional left and right factor of 
some iterate of $A$. The following result provides a partial converse statement (\cite{tame}).

\bt \l{cor} Let $A$ and $B$ be tame rational functions. 
Then the following conditions are equivalent.

\begin{enumerate} [label=\arabic*)]

\item Each iterate of $B$ is a compositional left  factor of 
some iterate of $A$. 

\item Each iterate of $B$ is a compositional right  factor of 
some iterate of $A$.

\item The functions $A$ and $B$ have a common iterate.  \qed 
\end{enumerate} 
\et

For rational functions $A$ and $B$, let us define an algebraic curve $  \f C_{A,B}$ by the formula $$\f   C_{A,B}:\, A(x)-B(y)=0.$$ The last result about tame rational functions we need below is following (see \cite{tame}, Corollary 3.6). 

\bt \l{c2} Let $A$ and $B$ be rational functions such that the curve 
$\f C_{A^{\circ s}, B},$ $s\geq 1,$ has an irreducible factor $\f C$ of genus zero or one for some $s\geq 1$. 
Assume in addition that $B$ is tame, $\deg A\geq 2,$  and  
\be \l{vod} s>\log_2\big[84(\deg B-1)(\deg B)!\big]. \ee
 Then  $A^{\circ s}=B\circ Q$ for some rational function $Q$, and 
$\f C$ is the graph $Q(x)-y=0.$ \qed 
\et

\subsection{Stabilization of semigroups $C(P^{\circ s})$}
For a rational function  $P$ of degree at least two, the 
groups in the sequence $ G(P^{\circ k})$, $k\geq 1,$ in general are different. 
Nevertheless, the following statement holds (\cite{sym}).

\bt \l{maa} Let $P$ be a rational function of degree $n\geq 2$.
Then  the sequence  $ G(P^{\circ k})$, $k\geq 1,$ contains only finitely many non-isomorphic groups, and, unless $P$ is a quasi-power, the orders of these groups    
are finite and uniformly bounded  in terms of   $n$ only. \qed
\et

Among other things, Theorem \ref{maa} implies that, unless $P$ is conjugate to $z^{\pm n},$ the group  $\Aut_{\infty}(P)$ is finite, so that  
\be \l{gopik} \Aut_{\infty}(P)=\Aut(P^{\circ s})\ee for some $s\geq 1$ (see \cite{sym} for more detail).  
In this section, we prove the following generalization of  equality \eqref{gopik}  for tame rational functions.

\bt \l{hh} Let $P$ be a tame rational function. 
Then $C_{\infty}(P)=C(P^{\circ s})$ for some $s\geq 1.$ 
\et 

\pr 
 Assume that $F\in C_{\infty}(P)$. Then, by the Ritt theorem, $F$ is a compositional right factor of some iterate of $P$. On the other hand, by Theorem \ref{1+}, there exist rational functions $F_1,F_2,\dots, F_s$ such that any compositional right factor of an iterate of $P$ has the form \eqref{egik}.  Furthermore, by Lemma \ref{per}, the function $\alpha\circ F_i\circ P^{\circ l}$  commutes with $P^{\circ s},$ $s\geq 1,$ if and only if $\alpha\circ F_i$ commutes with $P^{\circ s}$.

Let us observe now that if $\alpha\circ F_i$ commutes with $P^{\circ s}$, and $\alpha'\circ F_i$ commutes with $P^{\circ s'}$ for some $\alpha,\alpha'\in \Aut(\C\P^1)$ and $s,s'\geq 1,$  then the both functions $\alpha\circ F_i$ and $\alpha'\circ F_i$ commute with $P^{\circ \LCM(s,s')}$. 
Therefore, since 
$$\alpha'\circ F_i=(\alpha'\circ \alpha^{-1})\circ \alpha\circ F_i,$$
 Lemma \ref{per} implies that  $\alpha'\circ \alpha^{-1}$ also  commutes with $P^{\circ \LCM(s,s')}$. Thus, $\alpha'=\nu\circ \alpha$ for some $\nu \in \Aut_{\infty}(P)$. 
Since the group $\Aut_{\infty}(P)$ is finite, this yields that 
there exist finitely many rational functions $G_1,G_2,\dots, G_r\in C_{\infty}(P)$ such that  
 $F$ belongs to $C_{\infty}(P)$ if and only if $F$ has the form $$G_i\circ P^{\circ l}, \ \ \ 1\leq i \leq r, \ \ \ l\geq 0.$$
Finally, if $G_i$, $1\leq i \leq r,$ commutes with $P^{\circ k_i}$, $k_i\geq 1,$ then $G_i$ also commutes with $P^{N},$ where $N=\LCM(k_1,k_2,\dots,k_r).$  
Thus, $C_{\infty}(P)\subseteq C(P^{\circ N}),$ implying that 
$C_{\infty}(P)=C(P^{\circ N}).$ \qed

\subsection{Amenable semigroups} 
The following result is an extended version of Theorem \ref{11} from the introduction.

\bt \l{111} Let $S$ be a tame semigroup of rational functions. Then the following conditions are equivalent.

\begin{enumerate} [label=\arabic*)]

\item The semigroup $S$   is    left reversible. 

\item The semigroup $S$   is   left  amenable.

\item The semigroup $S$   is   amenable.

\item  The semigroup $S$   is  a subsemigroup of $C(P)$  
  for some tame rational \linebreak function $P$.

\item The semigroup $S$   is  power joined. 

\item The semigroup $S$ is left or right Archimedean.

\item The semigroup $S$ contains no free subsemigroup of rank two.

\item For all $A,B\in S$ there exist $z_1,z_2\in \C\P^1$ such that the forward orbits $O_A(z_1)$ and 	 $ O_B(z_2)$ have an infinite intersection.

\end{enumerate} 

\et

\pr
If $S$ is a power joined semigroup of rational functions, then   $S$ is a subsemigroup of 
 $C_{\infty}(F)$ for every $F\in S,$ by Theorem \ref{3}. Moreover, if $S$ is tame, then  every 
 $F\in S$ has degree at least two and is not special, since special rational functions are wild (see \cite{tame}, Corollary 2.5).  Therefore, by Theorem \ref{hh},   $C_{\infty}(F)=C(P)$, where $P=F^{\circ s}$ for some $s\geq 1.$ 
 Since $F^{\circ s}$ is tame by Lemma \ref{ltam}, this proves the 
 implication $5\Rightarrow 4$.  The implication $4\Rightarrow 3$ holds by Theorem \ref{2}. 
 The implication $3\Rightarrow 2$ is clear. The implication $2\Rightarrow 1$ holds by Proposition \ref{pater}. 

The  implication $1\Rightarrow 5$ follows from Theorem \ref{cor} and Theorem \ref{c2}. Indeed, let $A$ and $B$ be arbitrary elements of $S$. It follows from the left reversibility of $S$ that for every $s\geq 1$ there exist $C_{s},D_{s}\in S$ such that the equality 
$A^{\circ s}\circ  C_{s}= B\circ D_{s}$ 
holds, implying that the curve 
$\f C_{A^{\circ s}, B}$ has an irreducible factor of genus zero.    Since for $s$ big enough inequality \eqref{vod} holds, it follows from Theorem \ref{c2} that the function $B$  
is a compositional left  factor of some iterate of $A$  (notice that this fact does not immediately imply that $S$ is right Archimedean, since in the 
equality $A^{\circ n}=B\circ X$ the function $X$ may not belong to $S$). Moreover, 
using the same reasoning for iterates of $B$ we conclude that each iterate of $B$ is a compositional left  factor of 
some iterate of $A$, implying that $A$ and $B$ have a common iterate by Theorem \ref{cor}.  
This finishes the proof of the equivalences $1\Leftrightarrow 2\Leftrightarrow 3\Leftrightarrow 4\Leftrightarrow 5$.

The implication $5\Rightarrow 6$ follows from Lemma \ref{svin1}. 
On the other hand, if $S$ is left (resp. right) Archimedean, then for all $A$, $B\in S$ each iterate of $B$ is a compositional right (resp. left) factor of some iterate of $A$, implying by Theorem \ref{cor} that $A$ and $B$ share an iterate. Thus, $5\Leftrightarrow 6$.       
Further, it is clear that $5\Rightarrow 8$. On the other hand, it was proved in \cite{tame} that if  for 
tame rational functions $A$ and $B$ there exist orbits $O_A(z_1)$ and $ O_B(z_2)$ with an infinite intersection, then $A$ and $B$ share an iterate. For the reader convenience, we repeat the proof which relies on the Faltings theorem combined with Theorem \ref{cor} and Theorem \ref{c2}. We recall that by the Faltings theorem (\cite{fa}) if an irreducible algebraic curve $C$ defined over a finitely generated field $K$ of characteristic zero has infinitely many $K$-points, 
then $g(C)\leq 1.$  On the other hand, it is easy to see that if $O_A(z_1)\cap O_B(z_2)$ is infinite,  then for every pair $(i,j)\in \N\times \N$ the algebraic curve 
\be \l{xer} 
A^{\circ i}(x)-B^{\circ j}(y)=0\ee has infinitely many points $(x,y)\in O_A(z_1)\times O_B(z_2).$ 
Defining now $K$ as the field generated over $\Q$ by $z_1$, $z_2$, and the coefficients of $A$, $B$, and observing that the orbits $O_A(z_1)$ and $ O_B(z_2)$ belong to $K$, we conclude 
that for every pair $(i,j)\in \N\times \N$     curve \eqref{xer} has 
a factor of genus zero or one. It follows now from Theorem \ref{c2} that 
each iterate of $B$ is a compositional left  factor of some iterate of $A$, implying that $A$ and $B$ have a common iterate by Theorem \ref{cor}. Thus, $5\Leftrightarrow 8$

Since any tame semigroup of rational functions is left cancellative,  the implication $7\Rightarrow 1$ follows from Lemma \ref{xrun}. Finally, let us observer that a subsemigroup $S$ of $C(P)$ cannot contain a free subsemigroup of rank two $S'$, since such $S'$ is also a subsemigroup of $C(P)$ and hence is power joined by Theorem \ref{3}.   Thus, $4\Rightarrow 7$.   \qed 

\section{Archimedean and reversible semigroups of polynomials} 
\subsection{Functional equations in polynomials} 
We recall that a polynomial $A$ is called special if it is conjugate to $z^n$, $n\geq 2$  or to $\pm T_n$, $n\geq 2.$ The following two lemmas follow easily from from the characterization of the polynomials $z^n$ and $T_n$ in terms of their ramification  (see e.g. \cite{gtz2}, Lemma 3.5 and Lemma 3.9). 

\bl \l{any} Any decomposition of $z^n$, $n\geq 2,$ into a composition of polynomials  has the form \be \l{ega} z^n=(z^{n/d}\circ \mu)\circ(\mu^{-1}\circ  z^d),\ee where $d\vert n$ and $\mu$ is a polynomial of degree one. On the other hand, 
 any decomposition of $T_n$, $n\geq 2,$  has the form \be \l{xerxyi} T_n=(T_{n/d}\circ \mu)\circ(\mu^{-1}\circ  T_d),\ee where $d\vert n$ and $\mu$ is a polynomial of degree one. \qed
\el

For brevity, we will say that two polynomials $A$ and $B$ are {\it linearly equivalent} if there exist 
polynomials of degree one $\sigma$ and $\nu$ such that the equality $$A=\sigma\circ B\circ \nu$$ 
holds.

\bl \l{lemm} 
Let $A$ be a polynomials of degree $d\geq 2$ such that $A^{\circ l}$, $l> 1$, is linearly equivalent to $z^{d^l}.$ Then $A$ is conjugate to $z^d$. Similarly, if  $A^{\circ l}$, $l> 2$, is linearly equivalent to $T_{d^l},$ then $A$ is conjugate to $\pm T_d$. \qed 
\el

Since a rational function $A$ is a polynomial if and only if $A^{-1}\{\infty\}=\infty$,  for any 
decomposition $A=U\circ V$  of a polynomial $A$ into a composition of rational functions $U$ and $V$, there exists a M\"obius transformation $\mu$ such that $U\circ \mu$ and $\mu^{-1}\circ V$ are polynomials. Thus, considering decompositions of polynomials into compositions of rational functions, we can restrict ourselves by the consideration of decompositions into compositions of polynomials. Furthermore, it is not hard to prove that if an iterate of a rational function $F$ is a polynomial, then either $F$ itself is a polynomial or 
$F$ is conjugate to $1/z^d.$ This implies, in particular, that if $P$ is a polynomial of degree at least two not conjugate to $z^n$, then the semigroup $C(P)$  consists of polynomials. 

\vskip 0.2cm

Unlike the general case, {\it polynomial}  solutions of the functional equation \be \l{ura} A\circ C=B\circ D\ee admit essentially a complete description. 

Specifically, the following result follows easily from the fact that the monodromy group of a polynomial of degree $n$ contains a cycle of length $n.$

\bt [\cite{en}]\l{r1}
Let $A,C,B,D$ be polynomials 
such that $A\circ C=B\circ D$. Then there exist polynomials
$U, V, \widetilde A, \widetilde C, \widetilde B, \widetilde D,$ where
$$\deg U=\GCD(\deg A,\deg B),  \ \ \ \deg V=\GCD(\deg C,\deg D),$$
such that
$$A=U\circ \widetilde A, \ \  B=U\circ \widetilde B, \ \ C=\widetilde C\circ V, \ \  D=\widetilde D\circ V,$$
and 
$$ \widetilde A\circ \widetilde C=\widetilde B\circ \widetilde D.\eqno{\Box}$$
\et

Notice that Theorem \ref{r1} implies that if $\deg B\mid \deg A$ in \eqref{ura}, then  the equalities 
$$A=B\circ R, \ \ \ \ D=R\circ C$$ hold for some polynomial $R$.

Theorem \ref{r1} reduces describing solutions of \eqref{ura} to  describing solutions satisfying 
\be \l{ggccdd} {\rm GCD}(\deg A, \deg B)=1, \ \ \ {\rm GCD}(\deg C,\deg D)=1.\ee
The following result called ``the second Ritt theorem" (\cite{r1}) describes such solutions.

\bt [\cite{r1}]\l{r2}
Let $A,C,B,D$ be polynomials  
such that \eqref{ura} and \eqref{ggccdd} hold.
Then there exist polynomials $\sigma_1,\sigma_2,\mu, \nu$ of degree one  
such that, up to a possible replacement of $A$ by $B$ and of $C$ by $D$, either
\begin{align}  \l{rp11}  &A=\nu \circ z^sR^n(z) \circ \sigma_1^{-1}, & &C=
\sigma_1 \circ z^n \circ \mu \\ 
&\l{rp11+} B=\nu \circ z^n \circ \sigma_2^{-1},& 
&D=\sigma_2 \circ
z^sR(z^n) \circ \mu, \end{align} 
where $R$ is a polynomial, $n\geq 1,$ $s\geq 0$, and $\GCD(s,n)=1,$ or
\begin{align} \l{rp21} &A=\nu \circ T_m \circ \sigma_1^{-1}, &
&C=\sigma_1 \circ T_n \circ \mu, \\
&\l{rp21+} B=\nu \circ T_n \circ \sigma_2^{-1}& &D=\sigma_2 \circ
T_m\circ \mu, \end{align} 
where $T_n, T_m$ are the Chebyshev polynomials, $n,m\geq 1$, and $\GCD(n,m)=1.$ 
\et

Theorem \ref{r2} implies the following corollary. 

\bc \l{kal} Let $A,C,B,D$ be polynomials  
such that \eqref{ura} and \eqref{ggccdd} hold and $\deg A>\deg B$. Then $B$ and $C$ 
are linearly equivalent either to powers or to Chebyshev polynomials. \qed \ec


\subsection{Decompositions of iterates} 
Below we will use  the polynomial versions of Theorem \ref{1+}, Theorem \ref{cor}, and Theorem \ref{hh} given below. These versions are more precise since they hold for all non-special polynomials, not only for tame ones.

The more precise version of Theorem \ref{1+} for polynomials is the following result  (see \cite{mz}, and also \cite{pj}, \cite{tame}).

\bt \l{1+1} Let $A$ be a polynomial of degree $n\geq 2$ not conjugate to $z^n$ or $\pm T_n.$ 
Then there exists an integer $N$, depending on  $n$ only, such that  any decomposition of  $A^{\circ d}$ with $d\geq N$ is induced by a decomposition of $A^{\circ N}$. \qed 
\et

Theorem \ref{1+1} implies the following useful criterion. 

\bt \l{xur} Let $B$ be a polynomial  of degree at least two. Assume that there exists a sequence of polynomials $F_i,$ $i\geq 1$, such that: 

\begin{enumerate}  [label=\arabic*)]
\item  Each $F_i,$ $i\geq 1$, is a compositional left factor of some iterate of $B$. 
\item  Each $F_i,$ $i\geq 1$, is linearly equivalent to a special polynomial. 
\item  $\overline{\lim}_{i \to \infty}\deg F_i= \infty$. 
\end{enumerate}
Then $B$ is special. Moreover, the same conclusion holds if to replace the first condition by the condition that each $F_i,$ $i\geq 1$, is a compositional right factor of some iterate of $B$.
\et 
\pr We consider the ``left'' case. The proof in the ``right'' case is similar. Assume that $B$ is not special. 
Then Theorem \ref{1+1} implies that there exist a left compositional factor $C$ of some iterate of $B$ and different $i_1,i_2\geq 1$ such that 
\be \l{ert} F_{i_1}=B^{\circ l_1}\circ C\circ \mu_1, \ \ \ \ \ \ F_{i_2}=B^{\circ l_2}\circ C\circ \mu_2\ee for some 
$l_1,l_2\geq 1$ and $\mu_1,\mu_2\in \Aut(\C).$ Moreover, we can find $i_1$ and $i_2$ such that $l_1-l_2>2$. Equalities \eqref{ert} yield that 
$$ F_{i_1}=B^{\circ (l_1-l_2)}\circ F_{i_2}\circ (\mu_2^{-1}\circ \mu_1),$$ implying  by Lemma \ref{any} that $B^{\circ (l_1-l_2)}$ is linearly equivalent to a special polynomial. However, in this case $B$ is  special
by Lemma \ref{lemm}, in  contradiction with the assumption. \qed 

\vskip 0.2cm

The following analogue of Theorem \ref{hh} is obtained from Theorem \ref{1+1} in the same way as  Theorem \ref{hh} is obtained from Theorem \ref{1+}.

\bt \l{hh1} Let $P$ be a non-special polynomial of degree at least two. 
Then $C_{\infty}(P)=C(P^{\circ s})$ for some $s\geq 1.$ \qed 
\et 


Finally, the next two results are the ``left'' and the ``right'' polynomial analogues of Theorem \ref{cor}. The first of them was established previously  in the papers \cite{gtz}, \cite{gtz2} (see  \cite{gtz}, Proposition 3.3 and \cite{gtz2}, Proposition 4.1). The proof given below is somewhat shorter and is easily modified to fit the ``right'' case. It mimics the proof of Theorem 1.4 in \cite{tame}.

\bt \l{cor1} Let $A$ and $B$ be polynomials of degree at least two such that 
 each iterate of $B$ is a compositional left  factor of 
some iterate of $A$. Then either both $A$ and $B$ are special, or  there exist $k,l\geq 1$ such that $A^{\circ k}=B^{\circ l}$.
 \et
\pr Let us observe first that if one of the polynomials $A$, $B$ is special, then the other one is also special. Indeed, by condition,
for every $i\geq 1$ there exist $s_i\geq 1$ and $R_{i}\in \C[z]$ such that the equality
\be \l{zina} A^{\circ s_i}=B^{\circ i}\circ R_{i}\ee holds.  
Therefore, if $B$ is special, then  considering the sequence of special polynomials $F_i=B^{\circ i}$, $i\geq 1$, and applying Theorem \ref{xur}, we conclude that $A$ is special. On the other hand, if $A$ is a special, then each $B_i,$ $i\geq 1$, 
is linearly equivalent to a special polynomial by Lemma \ref{any}, implying that $B$ is special  by Lemma \ref{lemm}. 

Let us assume now that both $A$ and $B$ are not special. In this case, without loss of generality we may assume that the group $G(B)$ is finite. Indeed, 
it is easy to see that if $B$ is a quasi power, then $B^{\circ 2}$ is not  a quasi power, unless $B$ is conjugate to a power. Therefore, 
considering instead of $B$ its second iterate we may assume that $B$ is not a quasi-power, implying  by Theorem \ref{prim} that  $G(B)$ is finite. 
Furthermore, Lemma \ref{lemm} implies that all the groups $G(B^{\circ i}),$ $i\geq 1,$ are also finite. 

It follows from \eqref{zina} by Theorem \ref{1+1} that there exist a rational function $U$ and increasing sequences of non-negative integers $f_k$, $k\geq 0,$ and  $v_k$, $k\geq 0,$  such that 
\be \l{svin-} B^{\circ f_k}=A^{\circ v_{k}}\circ U\circ \eta_{k}, \ \ \ \ \ k\geq 0,\ee for some   
$\eta_k\in Aut(\C\P^1)$. In turn, this implies that 
there exists an increasing sequence of non-negative integers $r_k$, $k\geq 1,$ such that 
\be \l{svin} B^{\circ f_k}=A^{\circ r_{k}}\circ B^{\circ f_0}\circ \m_{k}, \ \ \ \ \ k\geq 1,\ee for some   
$\mu_k\in Aut(\C\P^1)$. Furthermore, since \eqref{svin} implies that for every $ k\geq 1$ the function  $B^{\circ f_0}\circ \m_{k}$ is a compositional right factor  of an iterate of $B$, 
there exist a rational function $V$ and 
an increasing sequence of non-negative integers   ${k_l},$  $l\geq 0,$  
such that 
$$B^{\circ f_0}\circ \m_{{{k_l}}}=\theta_l\circ V,\ \ \ \ \ l\geq 0,$$
for some $\theta_l\in Aut(\C\P^1),$ implying that 
$$B^{\circ f_0}\circ \m_{{{k_l}}}=\delta_l\circ  B^{\circ f_0}\circ \m_{{k_0}},\ \ \ \ \ l\geq 1,$$
for some $\delta_l\in Aut(\C\P^1).$

Clearly, the M\"obius transformations $ \m_{{k_l}}\circ \m_{{k_0}}^{-1},$ $l\geq 1,$ belong to the group $ G(B^{\circ f_0})$.  Therefore, the finiteness of $G(B^{\circ f_0})$ yields that      
$$\m_{{k_{l_2}}}\circ \m_{{k_0}}^{-1}=\m_{{k_{l_1}}}\circ \m_{{k_0}}^{-1}$$ for some $l_2>l_1$, 
implying that $\m_{{k_{l_2}}}=\m_{{k_{l_1}}}.$  
It follows now from \eqref{svin} that  
$$ B^{\circ f_{k_{l_2}}}=A^{\circ  (r_{k_{l_2}}-r_{k_{l_1}})} \circ B^{\circ f_{k_{l_1}}},
$$
 implying that 
\be \l{kro} B^{\circ (f_{k_{l_2}}-f_{k_{l_1}})}= A^{\circ (r_{k_{l_2}}-r_{k_{l_1}})} .\ee
Since $l_2>l_1$ and the sequences $k_l$, $l\geq 1,$ and $f_k,$ $k\geq 1,$ are increasing, 
the inequality $f_{k_{l_2}}>f_{k_{l_1}}$ holds, and therefore $A$ and $B$ have a common iterate. \qed 

\vskip 0.2cm
Since equality \eqref{skun} implies equality \eqref{pak}, if rational functions $A$ and $B$ satisfy condition \eqref{piz0}, then  each iterate of $B$ is a compositional right  factor of 
some iterate of $A$. Correspondingly, the ``right'' counterpart of Theorem \ref{cor1} is the following statement.

\bt \l{cor12} Let $A$ and $B$ be polynomials of degree at least two such that 
 each iterate of $B$ is a compositional right  factor of 
some iterate of $A$. Then either both $A$ and $B$ are special, or 
 there exist $k,l\geq 1$ such that   $A^{\circ 2k}=A^{\circ k}\circ B^{\circ l}$ and $B^{\circ 2l}=B^{\circ l}\circ A^{\circ k}$. 
\et
\pr In the same way as in the proof of Theorem \ref{cor1}, we conclude first that it is enough to prove the theorem assuming that both $A$ and $B$ are non-special and the group $G(B)$ is finite.  
The rest of the proof is obtained by a modification of the proof of Theorem \ref{cor1} as follows.  
Assuming that for every $i\geq 1$ there exist $s_i\geq 1$ and $R_{i}\in \C[z]$ such that \be \l{zona} A^{\circ s_i}= R_{i}\circ B^{\circ i},\ee we conclude that 
 there exists  a sequence $f_k$, $k\geq 0,$ such that 
\be \l{svon} B^{\circ f_k}=\m_{k}\circ B^{\circ f_0}\circ A^{\circ r_{k}} , \ \ \ \ \ k\geq 1,\ee for some   
$\mu_k\in Aut(\C\P^1)$ and $r_k\geq 1$, implying that there exists an increasing sequence ${k_l},$  $l\geq 0,$  
such that 
$$\m_{{{k_l}}}\circ B^{\circ f_0} =\m_{{k_0}}\circ  B^{\circ f_0}\circ \delta_l,\ \ \ \ \ l\geq 1,$$
for some $\delta_l\in Aut(\C\P^1)$. In turn, this yields that for some $l_2>l_1$ the equalities $\delta_{l_2}=\delta_{l_1}$ and  $\m_{{k_{l_2}}}=\m_{{k_{l_1}}}$ hold, implying by   
 \eqref{svon} that 
\be \l{since}  B^{\circ f_{k_{l_2}}}=B^{\circ f_{k_{l_1}}} \circ A^{\circ  (r_{k_{l_2}}-r_{k_{l_1}})}. 
\ee

Since \eqref{since} implies that each iterate of $A$ is a compositional right  factor of 
some iterate of $B$, repeating the above reasoning we conclude that \eqref{piz0} holds   for some $k_1,k_2,l\geq 1$. Arguing now as in the proof of Lemma \ref{ptw}, we conclude that \eqref{piz0} and \eqref{since} imply that there exist $k,l\geq 1$ such that the equalities  $A^{\circ 2k}=A^{\circ k}\circ B^{\circ l}$ and $B^{\circ 2l}=B^{\circ l}\circ A^{\circ k}$ hold.  \qed

\subsection{Reversible semigroups of polynomials}  
In this section, we characterize left or right reversible semigroups of polynomials by studying the corresponding functional equations.  

The following result was proved in the paper \cite{gtz2} (see \cite{gtz2},  Pro\-po\-sition 6.3). As above, we give an independent proof which can be modified to fit the ``right'' case.

\bt \l{bear} Let $A$ and $B$ be polynomials of degree $n\geq 2$ and $m\geq 2$  
 respectively such that for any $i,j\geq 1$ there exist  polynomial $C_{i,j},D_{i,j}$ satisfying \be \l{e11} A^{\circ i}\circ  C_{i,j}= B^{\circ j}\circ D_{i,j}.\ee 
Then either both  $A$ and $B$ are special, or there  exist $k,l\geq 1$ such that $A^{\circ k}=B^{\circ l}$.
\et
\pr For an integer  $n\geq 2$, let us denote by $\f P(n)$ the set of prime factors of $n$. 
Assume first that  \be \l{inc1} \f P(\deg B)\subseteq\f P(\deg A).\ee In this case, for every $j\geq 1$ the number $\deg B^{\circ j}$ is a divisor of the number $\deg A^{\circ i}$ for $i$ big enough. Therefore, by Theorem \ref{r1} applied to equality \eqref{e11}, for every $j\geq 1$ the polynomial $B^{\circ j}$ is a compositional left  factor of 
some iterate of $A$, implying by Theorem \ref{cor1} that either both  $A$ and $B$ are special, or $A$ and $B$ share an iterate.  
By symmetry, the same conclusion holds if \be \l{inc2} \f P(\deg A)\subseteq\f P(\deg B).\ee
 
Assume now that neither of conditions \eqref{inc1}, \eqref{inc2} holds. In this case,  there exist  $p_1\in \f P(\deg A)$ such that $p_1\not\in \f P(\deg B)$, and $p_2\in \f P(\deg B)$ such that $p_2\not\in \f P(\deg A).$  Applying Theorem \ref{r1} to equality \eqref{e11}, we can find polynomials 
$U_{i,j}, V_{i,j}, \widetilde A_{i,j}, \widetilde C_{i,j}, \widetilde B_{i,j}, \widetilde D_{i,j},$ where
\be \l{xorek} \deg U_{i,j}=\GCD(\deg A^{\circ i},\deg B^{\circ j}),  \ \ \ \deg V_{i,j}=\GCD(\deg C_{i,j},\deg D_{i,j}),\ee
such that
\be \l{xorek2} A^{\circ i}=U_{i,j}\circ \widetilde A_{i,j}, \ \  B^{\circ j}=U_{i,j}\circ \widetilde B_{i,j}, \ \ C_{i,j}=\widetilde C_{i,j}\circ V_{i,j}, \ \  D=\widetilde D_{i,j}\circ V_{i,j},\ee and 
\be \l{pis} \widetilde A_{i,j}\circ \widetilde C_{i,j}=\widetilde B_{i,j}\circ \widetilde D_{i,j}.\ee
Moreover, 
$$\gcd(\deg \widetilde A_{i,j}, \deg \widetilde B_{i,j})=1, \ \ \  \ \ \gcd(\deg \widetilde C_{i,j}, \deg \widetilde D_{i,j})=1,$$ and 
\be \l{gaf} \deg \widetilde A_{i,j}\geq p_1^i, \ \ \ \ \  \deg \widetilde B_{i,j}\geq p_2^j.\ee

Since the second equality in \eqref{xorek2} implies that 
$$\deg  \widetilde B_{i,j} \leq \deg B^{\circ j},$$ the degree  $\deg  \widetilde B_{i,j}$ is  bounded for fixed $j$. On the other hand, the first inequality in \eqref{gaf} implies that  $\deg \widetilde A_{i,j}\to \infty$ as $i\to \infty$. 
Therefore, applying Corollary \ref{kal} to equality \eqref{pis} for fixed $j$ and $i=i(j)$ big enough, we see that $\widetilde B_{i,j}$ is linearly equivalent to a special polynomial. 
It follows now from the second inequality in \eqref{gaf} and the second equality in \eqref{xorek2} that there exists a sequence of  polynomials $F_j$, $j\geq 1$, where $F_j=\widetilde B_{i,j}$ for some $i=i(j)$, satisfying conditions of Theorem \ref{xur}. Thus, $B$ is special. Moreover, by symmetry, $A$ is also special.  \qed

\bt \l{bear2} Let $A$ and $B$ be polynomials of degree $n\geq 2$ and $m\geq 2$  
 respectively such that for any $i,j\geq 1$ there exist  polynomial $C_{i,j},D_{i,j}$ satisfying \be \l{agai} C_{i,j}\circ A^{\circ i} =D_{i,j}\circ  B^{\circ j}.\ee
Then either both  $A$ and $B$ are special, or  there exist $k,l\geq 1$ such that   \linebreak $A^{\circ 2k}=A^{\circ k}\circ B^{\circ l}$ and $B^{\circ 2l}=B^{\circ l}\circ A^{\circ k}$.
\et
\pr 
If at least one of conditions \eqref{inc1}, \eqref{inc2} holds, then 
modifying the proof of Theorem \ref{bear}  using Theorem \ref{r1} and Theorem \ref{cor12},   we conclude that the theorem is true. 

On the other hand, if 
neither of conditions \eqref{inc1}, \eqref{inc2} holds, then  we can find polynomials 
$U_{i,j}, V_{i,j},$  $\widetilde A_{i,j}, \widetilde C_{i,j},$ $ \widetilde B_{i,j}, \widetilde D_{i,j},$ where
\be \l{xorek+} \deg U_{i,j}=\GCD(\deg C_{i,j},\deg D_{i,j}),  \ \ \ \deg V_{i,j}=\GCD(\deg A^{\circ i},\deg B^{\circ j}),\ee
such that
\be \l{xorek2+} C_{i,j}=U_{i,j}\circ \widetilde C_{i,j}, \ \  D_{i,j}=U_{i,j}\circ \widetilde D_{i,j}, \ \ A^{\circ i}=\widetilde A_{i,j}\circ V_{i,j}, \ \  B^{\circ j}=\widetilde D_{i,j}\circ V_{i,j},\ee 
\be \l{pis+} \widetilde C_{i,j}\circ \widetilde A_{i,j}=\widetilde D_{i,j}\circ \widetilde B_{i,j},\ee
and inequalities \eqref{gaf} hold for some primes $p_1,p_2$. Now a modification of the  proof of Theorem \ref{bear} shows that $A$ and $B$ are special. \qed 

\vskip 0.2cm

Theorem \ref{bear} and Theorem \ref{bear2} imply the following characterizations of 
left and right reversible semigroups. 

\bc \l{to3} 

Let $S$ be a semigroup of  polynomials of degree at least two containing at least one non-special polynomial.  Then the following conditions are equivalent.

\begin{enumerate} [label=\arabic*)]

\item The semigroup $S$   is    left reversible.

\item The semigroup $S$ is right Archimedean.

\item The semigroup $S$   is  power joined.

\end{enumerate} 

\ec 
\pr 
 The implications $3\Rightarrow 2$ and $2\Rightarrow 1$  follow from Lemma \ref{svin1}.   Finally, in view of Theorem \ref{bear}, to prove the implication $1\Rightarrow 3,$ it is enough to show  that if 
a left reversible semigroup of  polynomials $S$ contains a non-special polynomial $P$, then all elements of $S$ are non-special. In turn, the last statement also follows from Theorem \ref{bear}
applied to an arbitrary element $A$ of $S$ and $B=P$. \qed

\bc \l{to4}

Let $S$ be a semigroup of  polynomials of degree at least two containing at least one non-special polynomial.  Then the following conditions are equivalent.

\begin{enumerate} [label=\arabic*)]

\item The semigroup $S$   is    right reversible.

\item The semigroup $S$ is left Archimedean.

\item The semigroup $S$   is  power twisted. 

\end{enumerate} 

\ec
\pr 
 The proof  is similar to the proof of Corollary \ref{to3} with the use of Theorem \ref{bear2} instead of Theorem \ref{bear}. \qed

\vskip 0.2cm

Notice that Corollary  \ref{to3} and Corollary  \ref{to4} are not true for semigroups containing special polynomials. Indeed, for example, 
the semigroup of all Chebyshev polynomials is commutative and therefore is left and right inversible. 
However,  not all Chebyshev polynomials share an iterate.    
Similarly, one can easily see that the semigroup $S$ generated by $T_6$ and $T_{12}$, say, is left and right Archimedean, but is not power joined or power twisted since $6^k\neq 12^l$ for all $k,l\geq 1.$  
The assumption that 
$S$ contains only polynomials of degree at least two is also essential, since as it was already mentioned above any power joined semigroup $S$ with non-empty $\underline{S}$ coincides with $\underline{S}$.

\section{Amenable semigroups of polynomials} 

\subsection{Semigroups $C(P)$}  In this section, we describe  in terms of semidirect products the structure of the semigroup $C(P)$ for a non-special polynomial $P$. We will deduce this description from the following result  (see \cite{pj}, Theorem 1.3).

\bt \l{uni}
Let $P$ and $B$ be fixed non-special polynomials of degree at least two, and  let $\f E(P,B)$ be the set of all polynomials of degree at least two $X$ such that $P\circ X=X\circ B$.
Then, either $\f E(P,B)$ is empty, or there exists $R\in \f E(P,B)$ such that a polynomial $X$ belongs to 
$\f E(P,B)$ if and only if $X= A\circ R$ for some polynomial  $A$ commuting with $P.$ \qed
\et

Theorem \ref{uni} implies the following statement.

\bt \l{co4} Let $P$ be a non-special polynomial of degree at least two. Then for every $A\in \overline{C(P)}$ 
 the group $Aut(P)$ is a subgroup 
of the group $G(A)$ such that  $\gamma_A(\Aut(P))= \Aut(P)$. Furthermore, there exists a polynomial $R\in \overline{C(P)}$  such that  $C(P)=S_{\Aut(P),R}.$

\et 
\pr 
Since $\f E(P,P)=C(P)$, it follows from Theorem \ref{uni} that  there exists \linebreak $R\in C(P)$ such that every $A\in C(P)$  has the form 
\be \l{inn} A=U\circ R,\ee where $U\in C(P)$. In turn, this implies that every $A\in C(P)$  can be represented in the form  
\be \l{form} A=\sigma \circ R^{\circ s}\ee for some $s\geq 0$ and $\sigma\in \Aut(P)$. Indeed, if $\deg U=1$ in \eqref{inn}, then \eqref{form} holds for $s=1$. Otherwise, applying Theorem \ref{uni} again to the polynomial $U$ and so on, we obtain \eqref{form} for some $s>1$.  

Further, if $A\in \overline{C(P)}$ and $\omega \in \Aut(P)$, then 
 $A\circ \omega\in \overline{C(P)}$, and using the representation \eqref{form} we see that for every $\omega \in \Aut(P)$ there exists 
$\delta \in \Aut(P)$ such that 
$$\delta\circ A=A\circ \omega.$$
Moreover, since $C(P)$ is cancellative by Theorem \ref{1}, to different $\omega$ correspond different $\delta$. Thus, for every $A\in C(P)$ the group $\Aut(P)$ is a subgroup 
of $G(A)$ such that $\gamma_A(\Aut(P))= \Aut(P)$. Finally, it follows from \eqref{form} that 
$C(P)\subseteq S_{\Aut(P),R},$ and it is clear that $S_{\Aut(P),R}\subseteq C(P)$. Hence,
 $S_{\Aut(P),R}= C(P)$. \qed

\bc \l{suka} Let $P$ be a non-special polynomial of degree at least two, and $A_1,A_2\in C(P)$ polynomials such that $\deg A_2\geq \deg A_1$. Then there exists a uniquely defined polynomial $U$ such that 
$A_2=U\circ A_1.$ Moreover, $U\in C(P)$. 
\ec 
\pr Indeed, 
since $$A_1=\sigma_1 \circ R^{\circ s_1}, \ \ \ \ A_2=\sigma_2 \circ R^{\circ s_2},$$
for some $s_1,s_2\geq 0$ and $\sigma_1,\sigma_2\in \Aut(P)$ the equality 
$A_2=U\circ A_1$ holds for $$U=\sigma_2\circ R^{\circ (s_2-s_1)} \circ \sigma_1^{-1}\in C(P).$$ The uniqueness of $U$ follows from the right cancellativity of $\C[z].$ \qed 

\vskip 0.2cm
Let us recall that the classical theorem of Ritt about commuting  polynomials (see \cite{rr}, \cite{r}) states that if $P_1$, $P_2$ is a pair of commuting polynomials of degrees $d_1\geq 2$ and $d_2\geq 2$, then up to the change 
\be \l{chan} P_1 \rightarrow \lambda \circ P_1 \circ \lambda^{-1}, \ \ \ \ P_2 \rightarrow \lambda \circ P_2\circ \lambda^{-1},\ee where $\lambda$ is a polynomial of degree one, either  
$$P_1=z^{d_1}, \ \ \  P_2=\v z^{d_2},$$ 
where $\v$ is a $(d_1-1)$th root of unity, or 
$$P_1=\pm T_{d_1}, \ \ \  P_2=\pm T_{d_2},$$ 
where $T_{d_1}$ and $T_{d_2}$ are the Chebyshev polynomials, or  
\be \l{thco} P_1=\v_1R^{\circ s_1},\ \ \  P_2=\v_2R^{\circ s_2},\ee
where $R=zS(z^{\ell})$ for some polynomial $S$ and $l\geq 1$, and $\v_1,\v_2$ are $l$th root of unity.

Notice that condition \eqref{thco} is equivalent to the condition that there exists a polynomial $R$ such that 
\be \l{rfv} P_1=\sigma_1\circ R^{\circ s_1},\ \ \  P_2=\sigma_2\circ R^{\circ s_2},\ee where 
$\sigma_1,\sigma_2\in \Aut(R).$ Indeed,  every polynomial is conjugate to a polynomial of the form \be \l{exx} R=z^n+c_{n-2}z^{n-2}+\dots + c_0,\ee where $c_n=1$ and $c_{n-1}=0$. Furthermore, one can easily see that if  $$\sigma=az+b$$ commutes with such $R$, then $b=0$ 
and $a$ is a root of unity, implying that $R=zS(z^{\ell})$ for some polynomial $S$ and $l\geq 1$, and $a$ is an $l$th root of unity. 

Theorem \ref{co4} implies  the Ritt theorem about commuting  polynomials. Moreover, it implies the classification of commutative semigroups of polynomials obtained by Eigenthaler and  Woracek (\cite{wora}). To formulate the corresponding result, we introduce some notation. 
We denote by $\f Z$ the semigroup consisting of  polynomials  of 
the form $az^n,$ where $a\in \C^*$ and $n\geq 1,$ and by  $\f T$ the semigroup consisting of  polynomials of the form $\pm T_n,$ $n\geq 1.$  We say that two semigroups of polynomials $S_1$ and $S_2$ are {\it conjugate} if 
there exists $\alpha\in \Aut(\C)$ such that $\alpha\circ S_1\circ \alpha^{-1}=S_2.$ 

\bt \l{coma1}  Let $S$ be a commutative semigroup of  polynomials containing at least one non-special polynomial of degree greater than one. Then either $S$ is 
conjugate to a subsemigroup of $\f Z$ or  $\f T$, or $S$  is 
a subsemigroup of $S_{\Aut(R),R}$ for some non-special polynomial $R$  of degree at least two. 
\et 
\pr
In case $S$ contains at least one special polynomial $P$, the theorem follows from the following simple fact: if the unite circle is a completely invariant set for a polynomial $P$, then   $P=az^n$, where $\vert a\vert =1$  (see \cite{be3}, Theorem 1.3.1). Since commuting polynomials have the same Julia sets, this implies that if some $P\in S$ is conjugate to $z^{n}$, then $S$ is 
conjugate to a subsemigroup of $\f Z$. 
Similarly, since $\pm T_n$ are the only polynomials whose Julia set is the unit segment,  if  some $P\in S$  is conjugate to $\pm T_{n}$, then $S$ is 
conjugate to a subsemigroup of $\f T$.

The above shows that we can assume that $S$ contains no special polynomials.  Let us consider the semigroup 
$$\f F=\bigcap_{P\in S}C(P).$$
 Since the semigroup $S$ belongs to $\f F$ by condition, $\f F$ contains polynomials of degree greater than one, and therefore contains a polynomial of minimum possible degree greater than one. Let $R$ be any such a polynomial, and $P$ an arbitrary element of $S$.   
Since $P,R\in \f F$ and $\deg P\geq \deg R$ by construction,  Corollary \ref{suka} yields that there exists a polynomial $Q\in \f F$ such that $P=Q\circ R$. Moreover, since the inequality $1<\deg Q< \deg R$ contradicts to the choice of $R$, one of the equalities 
 $\deg Q=1$ or $\deg Q\geq \deg R$ holds. In case $\deg Q\geq \deg R,$ we can apply the same reasoning to $Q$ and so on, eventually obtaining a representation 
$P=\mu\circ R^{\circ s}$ for some  $s\geq 1$ and a polynomial of degree one $\mu\in \f F$.  Furthermore, since $P,R\in C(R)$, it follows from Lemma \ref{per} that $\mu\in C(R).$ Thus, every $P\in S$ belongs to $S_{\Aut(R),R}$, and hence $S$ is a subsemigroup of $S_{\Aut(R),R}$. \qed

\vskip 0.2cm
Notice 
that the fact that Theorem \ref{uni} implies the Ritt theorem was already mentioned in the paper \cite{pj}. However, the proof given there is not complete since it provides a representation \eqref{rfv},  where actually  only $\mu_2$ belongs to $\Aut(R)$, while $\mu_1$ belongs to $\Aut(P)$. To correct it, one has to define $R$ as a polynomial of minimum possible degree commuting with {\it both} $P_1$ and $P_2,$ as above.

\subsection{Semigroups $I(K)$}
For a compact set $K\subset \C$, we denote by $I(K)$ the set of polynomials $A$ such that $A^{-1}\{K\}=K,$ and by $\Omega_K$ the subset of $I(K)$ consisting of polynomials of degree one.  It is clear that  $I(K)$ is a semigroup, and $\Omega_K$ is a group.  
In this section, we describe in terms of semidirect products the structure of the semigroup $I(K)$ for a compact set $K$ that is neither a union of concentric circles nor a segment.  Our approach is based on the following result from the paper \cite{p1} (see \cite{p1}, Theorem 3, and also the related papers \cite{d}, \cite{d2}). 

\bt \l{prei} Let  $K\subset \C$ 
be a compact set that is neither  a union of concentric circles nor a segment, and 
 $A_1, A_2$ polynomials of degree greater than one such that \be \l{xx} A_1^{-1}\{K\}=A_2^{-1}\{K\}=K. \ee 
Then  the group $\Omega_K$ is finite and 
there exists a polynomial $F$  such that $F^{-1}\{K\}=K$ and
\be \l{po3}
A_1= \mu_1 \circ F^{\circ s_1} , \ \ \ A_2=\mu_2 \circ F^{\circ s_2}
\ee
for some $\mu_1, \mu_2\in\Omega_K$ and $s_1,s_2\geq 0.$   \qed 
\et 

Theorem \ref{prei} implies the following result. 

\bt \l{co3}  Let  $K\subset \C$ be a compact set that is neither a union of concentric circles nor a segment such that $\overline{I(K)}\neq \emptyset.$ 
Then for every $A\in \overline{I(K)}$ the group $\Omega_K$ is a subgroup 
of the group $G(A)$ such that $\gamma_A(\Omega_K)\subseteq \Omega_K$. Furthermore, there exists $R\in \overline{I(K)}$  such that $I(K)=S_{\Omega_K,R}$. 
\et
\pr Let $\mu$ be an arbitrary element of $\Omega_K$. Then $A\circ \mu\in \overline{I(K)}$ for every $A\in \overline{I(K)}$,  implying by Theorem \ref{prei} that there exist $s\geq 1$ and $F\in \overline{I(K)}$ such that 
$$A= \mu_1 \circ F^{\circ s} , \ \ \ A\circ \mu=\mu_2 \circ F^{\circ s}$$ for some $\mu_1, \mu_2\in\Omega_K$.  Therefore, $$ A\circ \mu=\delta\circ A$$
where $\delta  =\mu_2\circ \mu_1^{-1}\in \Omega_K$, and 
 hence $\gamma_A(\Omega_K)\subseteq \Omega_K$.

Further, it is clear that $S_{\Omega_K,R}\subseteq I(K)$ for every element $R$ of $ \overline{I(K)}$. 
 On the other hand, if $R$ is any polynomial of minimum possible degree which belong to $\overline{I(K)}$, then Theorem \ref{prei} implies that for every $A\in \overline{I(K)}$   there exists $F\in I(K)$ such that  the equalities 
$$A= \mu_1 \circ F^{\circ s} , \ \ \ R=\mu_2 \circ F$$ hold for some $\mu_1, \mu_2\in\Omega_K$ and $s\geq 1,$ implying that 
 $$A=\mu_1 \circ (\mu_2^{-1} \circ R)^{\circ s}\in S_{\Omega_K,R}.\eqno{\Box}$$ 

Theorem \ref{co3} can be regarded as a generalization of the classification of  pairs of polynomials  sharing Julia sets (see \cite{a2}, \cite{b}, \cite{be2}, \cite{be1}, \cite{sh}). In particular, since for $R$ of the form \eqref{exx} the group $G(R)$  is non-trivial if and only if  
 $R=z^rS(z^{\ell})$ for some polynomial $S$ and $l\geq 1$, $r\geq 0$, Theorem \ref{co3} implies that if $P_1$ and $P_2$ is a pair of non-special polynomials of degrees at least two sharing a Julia set, then, up to the change \eqref{chan}, 
\be \l{thcc} P_1=\v_1R^{\circ s_1},\ \ \  P_2=\v_2R^{\circ s_2},\ee
where $R=z^rS(z^{\ell})$ for some polynomial $S$ and $l\geq 1$, $r\geq 0$, and $\v_1,\v_2$ are $l$th root of unity.

\subsection{Left amenable semigroups} 

The following result is an extended version of Theorem \ref{-1} for semigroups of  polynomials of degree at least two. 

\bt \l{22} 
Let $S$ be a semigroup of  polynomials of degree at least two containing at least one non-special polynomial.  Then the following conditions are equivalent.

\begin{enumerate} [label=\arabic*)]

\item The semigroup $S$   is    left reversible. 

\item The semigroup $S$   is   left amenable.

\item The semigroup $S$   is   amenable.

\item The semigroup $S$ is a subsemigroup of $\overline{S_{\Gamma,R}}$ for some non-special polynomial $R$   
of degree at least two and a subgroup $\Gamma$ of $G(R)$ such that $\gamma_R(\Gamma)= \Gamma$.

\item  The semigroup $S$   is  a subsemigroup of $\overline{C(P)}$  
  for some non-special polynomial $P$ of degree at least two.

\item The semigroup $S$   is  power joined.

\item The semigroup $S$ is right Archimedean.

\item For all $A,B\in S$ there exist $z_1,z_2\in \C\P^1$ such that the forward orbits $O_A(z_1)$ and 	 $ O_B(z_2)$ have an infinite intersection.

\end{enumerate} 

\et 
\pr Let $Q$ be a non-special polynomial that belongs to $S$. It follows from Theorem \ref{3} that if $S$ is power joined, then $S$ is a subsemigroup of $\overline{C_{\infty}(Q)}$,  and Theo\-rem \ref{hh1} implies that $$\overline{C_{\infty}(Q)}=\overline{C(P)},$$ where $P=Q^{\circ s}$ for some $s\geq 1.$ Moreover, $P$ is non-special by Lemma \ref{lemm}. Thus, $6\Rightarrow 5$. The implications $5\Rightarrow 4$ follows from Theorem \ref{co4}. The implication $
4\Rightarrow 3$ follows from Theorem \ref{moh2}. The implications 
 $3\Rightarrow 2$ is obvious. The 
implication $2\Rightarrow 1$ follows from Proposition \ref{pater}. Since $1\Leftrightarrow 6 \Leftrightarrow 7$ by Corollary \ref{to3}, this shows that  the first seven conditions of the theorem are equivalent.  

 The implication $6\Rightarrow 8$ is clear. 
Finally, it was shown in the papers \cite{gtz}, \cite{gtz2} (see also \cite{tame} for another proof) that if $A$ and $B$ are polynomials of degree at least two having orbits with infinite intersection, then  $A$ and $B$  have a common iterate. 
Thus, $8\Rightarrow 6$.
\qed 

\vskip 0.2cm

\noindent{\it Proof of Theorem \ref{-1}.}  The implications 
$5\Rightarrow 4\Rightarrow 3\Rightarrow 2\Rightarrow 1$  are proved in the same way as the 
corresponding implications in Theorem \ref{22}. On the other hand,  by Lem\-ma \ref{zaqw0}, to prove the implication $1\Rightarrow 5$ it is enough to prove that  $\overline{S}$  is  a subsemigroup of $\overline{C(P)}$. Therefore, this implication is a corollary of the corresponding implication from Theorem \ref{22}. Finally, it is clear that the final statement of the theorem also follows from Theorem \ref{22}. 
  \qed

\subsection{Right amenable semigroups} 
The following result is the analogue of Theorem \ref{22} for right amenable semigroups. 
\bt \l{33} 
Let $S$ be a semigroup of  polynomials of degree at least two containing at least one non-special polynomial.  Then the following conditions are equivalent.

\begin{enumerate} [label=\arabic*)]

\item The semigroup $S$   is   right reversible. 

\item The semigroup $S$   is right amenable.

\item The semigroup $S$ is a subsemigroup of $\overline{S_{\Gamma,R}}$ for some non-special polynomial $R$  of degree at least two and a subgroup $\Gamma$ of $G(R)$ such that $\gamma_R(\Gamma)\subseteq \Gamma$.

\item The semigroup $S$ is a subsemigroup of  $\overline{I(K)}$ for some compact set $K\subset \C$, which is neither a union of concentric circles nor a segment.

\item The semigroup $S$ is a subsemigroup of  $\overline{E(P)}$ for some  non-special polynomial $R$   
of degree at least two.

\item  The semigroup $S$ is power twisted.

\item The semigroup $S$ is left Archimedean.

\item The semigroup $S$ contains no free subsemigroup of rank two.

\end{enumerate} 

\et
\pr  Let $P$ be a non-special polynomial that belongs to $S$. If  $S$ is power twisted, then $S$ is a subsemigroup of  $E(P)$ by Theorem \ref{kosh2}.  Therefore, $6\Rightarrow 5$.  Further, since for every rational function $P$ the support of $\mu_P$ coincides with $J(P)$, rational functions sharing the measure of maximal entropy with $P$ share the  Julia set with $P$, implying that   $\overline{E(P)}\subseteq \overline{I(J(P))}$. 
On the other hand, since the Julia set $J(P)$ of a polynomial $P$ is the boundary of its filled-in Julia set, if $J(P)$ is a union of circles, then $J(P)$  is a circle.   
Thus, taking into account  the above mentioned fact that the Julia set of a non-special polynomial $P$ cannot be a  circle or a segment, we conclude that  $5\Rightarrow 4$. The implication $4\Rightarrow 3$ follows from Theorem \ref{co3}. The implication $3\Rightarrow 2$ follows from Theorem \ref{moh1}. The 
implication $2\Rightarrow 1$ follows from Proposition \ref{pater}. 
Since $1\Leftrightarrow 6 \Leftrightarrow 7$ by Corollary \ref{to4}, this shows that  the first seven conditions of the theorem are equivalent.

The implication $8\Rightarrow 1$ follows from Lemma \ref{xrun}, since any semigroup of rational functions is right cancellative. Finally, to prove the implication  $3\Rightarrow 8$, we observe that if $S$ is a subsemigroup of $ S_{\Gamma,R}$, then every subsemigroup $S'$ of $S$ is also  a subsemigroup of $ S_{\Gamma,R}$. Therefore, $S'$  is right amenable by Theorem \ref{moh1} and hence   $S'$ is not  free by Theorem \ref{xrukas}. \qed

\vskip 0.2cm

\noindent{\it Proof of Theorem \ref{-11}.} 
In view of Theorem \ref{33},  to prove the implication $5\Rightarrow 4$ we only must show that $\underline{(E(P))}\subseteq I(J(P)).$  Let  $\mu$ be an element of $ \underline{E(P)}$.  Then, $P\circ \mu\in \overline{E(P)}$  and hence  $P\circ \mu\in I(J(P))$  by Theorem \ref{33}.  Thus, $J(P)=J(P\circ \mu)$,  and  the invariance of $J(P)$ with respect to $P$ and  $P\circ \mu$ implies  that 
$$J(P)=(P\circ \mu)^{-1}(J(P))= \mu^{-1}(P^{-1} (J(P))= \mu^{-1}(J(P)).$$ 
The implications 
$4\Rightarrow 3\Rightarrow 2\Rightarrow 1$ are proved in the same way as in the proof of Theorem \ref{33}. 
The implication $1\Rightarrow 5$ follows from Lemma \ref{zaqw}  combined with Theorem \ref{33}. The implications  $6\Rightarrow 1$ and $3\Rightarrow 6$ are obtained in the same way as the implications  $8\Rightarrow 1$ and $3\Rightarrow 8$  in Theorem \ref{33}. 
Finally,  the last statement of the theorem follows from Theorem \ref{33} and Lemma \ref{ptw}. 
\qed

\vskip 0.2cm
\noindent {\bf Acknowledgments}. The  author is grateful to Carlos Cabrera and Peter Makienko for interesting discussions.

\end{document}